\documentclass{amsart}

\newtheorem{thm}{Theorem}[section] 
\newtheorem{prop}[thm]{Proposition} 
\newtheorem{lem}[thm]{Lemma}

\begin{document}

\title[Irreducible free group representations]{Some irreducible free group representations in which a linear
combination of the generators has an eigenvalue}   
\author{William L. Paschke}
\address{Department of Mathematics \\
University of Kansas \\
Lawrence, KS 66045-2142 \\
USA}
\email{paschke@math.ukans.edu}

\begin{abstract} We construct irreducible unitary representations of a finitely generated free group which are
weakly contained in the left regular representation and in which a given linear combination of the generators has an
eigenvalue. When the eigenvalue is specified, we conjecture that there is only one such representation. The
representation we have found is described explicitly (modulo inversion of a certain rational map on euclidean space)
in terms of a positive definite function, and also by means of a quasi-invariant probability measure on the
combinatorial boundary of the group. 
\end{abstract}

\maketitle

\

\section{Statement of main result}\label{summary}

Let $G$ be the free group on $n$ generators $u_1, u_2, \ldots u_n,$ where $n \geq 2.$ For a reduced word
$s$ in $G$, let $|s|_i$ for each index $i$ (resp. $\gamma(s)$) be the number of occurrences of
$u_i$ or $u_i^{-1}$ (resp. $u_k^{-1} u_j$ for $k \not= j$) in $s$. Given positive numbers $c_1, \ldots, c_n$ and
a positive number $\lambda$ satisfying
$$ 2 \max c_i^2 - \sum_i c_i^2 < \lambda^2 < \sum_i c_i^2 ,$$
let $(x_1,  \ldots, x_n)$ be the unique vector in the positive orthant of ${\mathbb R}^n$ satisfying
$$\frac{x_j(1 + \sum_{i \not= j} x_i)}{\sum_i x_i} = \frac{c_j^2}{\lambda^2}$$
for each $j$. Define numbers $a_1, \ldots, a_n$ by
$$a_j = \sqrt{\frac{x_j}{\left(1 + \sum_{i \not= j} x_i \right) \sum_i x_i}}.$$
The function $\phi$ defined on $G$ by
$$\phi(s) = (- \sum_j x_j)^{\gamma(s)} \prod_j a_j^{|s|_j}$$
is positive definite and satisfies
$$\sum_i c_i \phi(s u_i) = \lambda \phi(s)$$
for all $s$ in $G$. The associated unitary representation $\pi$ is irreducible and weakly contained in
the left regular representation (so $\phi$ extends to a pure state of the reduced group $C^*$-algebra $C^*_r(G)$). The
$\lambda$-eigenspace of $\sum_i c_i \pi(u_i)$ is one-dimensional.

\

\section{Introduction}\label{intro} 

To make the discussion below more mellifluous, we will adopt the following terminology. By a state of $G$, we
mean a unital positive definite function on $G$. We will often regard a state  of $G$ as a positive linear
functional on the complex group algebra ${\mathbb C}G$, or on the universal $C^*$-algebra $C^*(G)$, or (if appropriate) on
the reduced $C^*$-algebra $C^*_r(G)$, that is, the $C^*$-algebra generated by the image of $G$ in its left regular
representation on $\ell^2(G)$. We call a state reduced if it is bounded with respect to the operator norm on
${\mathbb C}G$ in the left regular representation (and so extends to a state of $C^*_r(G)$). In general, we will use the
term reduced to mean ``having to do with the left regular representation,'' so for instance the reduced spectrum of an
element in ${\mathbb C}G$ is its spectrum in the left regular representation. We will call a unitary representation of $G$
reduced if it is weakly contained in the left regular representation, that is, if it extends to a $\ast$-representation of
$C^*_r(G).$ (We will avoid the locution reduced irreducible representation, however.) For a complex number
$\lambda$ and $X$ in ${\mathbb C}G$, a $\lambda$-eigenstate for $X$ is a state $\phi$ of $G$ such that $\phi((X^* -
\overline{\lambda})(X - \lambda)) = 0.$ This is equivalent to $\phi(s (X - \lambda)) = 0$ for all $s$ in $G$. If $\pi$ is
the unitary representation obtained from $\phi$ by the Gel'fand-Raikov construction and $\xi_0$ the corresponding
distinguished cyclic vector, then clearly $\phi$ is a $\lambda$-eigenstate for $X$ if and only if $\xi_0$ is a
$\lambda$-eigenvector for $\pi(X).$ Finally, we will call a state of $G$ pure if it is not a proper convex combination of
two different states of $G$. This is equivalent to irreducibility of the associated unitary representation.

Consider $Y$ in ${\mathbb C}G \setminus\{0\}$. If $Y$ is singular in $C^*_r(G)$, there is at least one reduced
$0$-eigenstate for $Y$; this is because unilateral and bilateral invertibility coincide in $C^*_r(G).$ We conjectured in
\cite{Paschke} that there are only finitely many reduced pure $0$-eigenstates for $Y$, and proved this when $Y =
u_1 + \ldots + u_n - \lambda$ with $|\lambda| = \sqrt{n}$ (in which case there is a unique reduced $0$-eigenstate for
$Y$),  and also when $Y$ is a polynomial in one of the generators (where the reduced pure eigenstates correspond to
distinct zeros of the polynomial on the unit circle). Kuhn and Steger \cite{KuhnSteger3}, using facts developed in
\cite{FigaTalamancaSteger}, have very recently shown that a selfadjoint linear combination $X$ of the generators and
their inverses has a unique reduced $\lambda$-eigenstate when $\lambda$ is plus or minus the reduced spectral
radius (= reduced norm) of $X$. The results in \cite{Paschke} and \cite{KuhnSteger3} resemble one another in
treating only extreme spectral values, but the techniques of proof are utterly different. 

At any rate, reduced pure eigenstates for linear combinations of group elements appear to be special enough that
there is some point to seeking them out in simple cases with a view to describing them and their associated
representations in detail.  In what follows, we consider a linear combination $X = c_1 u_1 + \ldots + c_n u_n$ of the
generators. The reduced spectrum of $X$ is the closed annulus with outer radius the euclidean length of the
coefficient vector $\vec{c}$, and inner radius either zero if no $|c_i|^2$ exceeds the sum of the other $|c_j|^2$'s, or else
the square root of the largest $|c_i|^2$ minus the sum of the other $|c_j|^2$'s. For $\lambda$ in the reduced spectrum,
the problem we address is mostly that of finding a reduced $\lambda$-eigenstate for $X$, but at the end of Section
\ref{further} we record another marginal uniqueness result. Namely, when $|\lambda| = |\vec{c}|$ the
argument in \cite{Paschke} goes over essentially {\em verbatim} to show that there is only one reduced
$\lambda$-eigenstate; the case in which $|\lambda|$ is the inner radius of the spectral annulus
(provided the latter is positive) then follows easily.  For other $\lambda$ in the reduced spectrum, we can only
conjecture uniqueness, but it is encouraging that the reduced $\lambda$-eigenstate we identify turns out to be pure,
and that the $\lambda$-eigenspace of $X$ in the associated representation turns out to be one-dimensional. 

In the case of equal coefficients treated in \cite{Paschke}, considerations of symmetry and economy give one a good
idea ahead of time of what the formula for a reduced $\lambda$-eigenstate should be.  The present case is more
suspenseful. We look for the desired needle in two different haystacks --- essentially, among $\lambda$-eigenstates
which are rarely reduced, and among reduced states which are rarely $\lambda$-eigenstates.  The states in the first
search venue are defined in terms of the functions $|\cdot|_i$ and $\gamma_{ij}$ that count, respectively the number
of $u_i$'s and $u_i^{-1}$'s in a reduced word, and the number of $u_i^{-1}u_j$'s and $u_j^{-1}u_i$'s. (Thus
$\gamma_{ij} = \gamma_{ji}$ and $\gamma_{ii} = 0$. Notice also that $\sum_i| \cdot |_i$ is the usual length function
on $G$ corresponding to the given generators.)  We assume henceforth that the coefficients $c_j$ are all positive (else
rotate by scalars of modulus one, and reduce the number of generators if any of the $c$'s vanish).  Write $\vec{c} =
(c_1, c_2, \ldots, c_n)$. Let $B = (b_{ij})$ be a positive $n \times n$ matrix  with real entries and with 1's on the
diagonal. Assume further (mostly for convenience) that $B \vec{c} \cdot \vec{c} > 0.$ Define the scalar
$\lambda$ and the vector $\vec{a} = (a_1, a_2, \ldots, a_n)$ by $\lambda = \sqrt{B \vec{c} \cdot \vec{c}}$ and
$\vec{a} = \lambda^{-1} B \vec{c}$. Define $\phi$  on $G$ by 
$$\phi(s) = \phi_{B,\vec{c}}(s) =  \prod_i a_i^{|s|_i - \sum_j \gamma_{ij}(s)} \prod_{i > j}
b_{ij}^{\gamma_{ij}(s)}$$ (where we read $0^0$ as 1). Thus, $\phi(u_j) = a_j, \phi(u_i u_j) = a_i a_j,
\phi(u_i u_j^{-1}) = a_i a_j \ ( i \not= j), \phi(u_i^{-1} u_j) = b_{ij} ( i \not= j), $ and so on. It turns out that
$\phi$ satisfies $\sum_i c_i \phi(s u_i) = \lambda
\phi(s)$ for all $s$ in $G$ (easy), is positive definite on $G$ (moderately difficult), and is pure
(strenuous), which is to say that the unitary representation of $G$ to which $\phi$ gives rise is
irreducible. (The proof of irreducibility uses a result of P. Linnell \cite{Linnell}, a special case of
which asserts that every nonzero element of the group algebra of $G$ has zero kernel in the left regular
representation. Linnell's clever proof of this involves the Fredholm module that is sometimes used in demonstrating
the absence of nontrivial idempotents in $C^*_r(G).$)

The $\lambda$'s that arise here are the positive numbers in the spectrum of $\sum c_i u_i$ in the
universal representation. The latter is the annulus with outer radius $\sum c_i$ and inner radius
$$\max\{0, 2 \max c_i - \sum c_i \}.$$ More restrictively, suppose that $\lambda$ is positive and lies in
the interior of the spectrum of $\sum c_i u_i$ in the left regular representation, that is,
$$ 2 \max c_i^2 - \sum c_i^2 < \lambda^2 < \sum c_i^2 .$$
There are many eligible matrices $B$ such that $B \vec{c} \cdot \vec{c} = \lambda^2$. Which one of these
--- and we imagine there can be at most one --- is such that $\phi_{B, \vec{c}}$ is reduced, that is,
extends to a state of $C^*_r(G)$?

At this point, the search moves to the second haystack, namely boundary representations of $G$. To
obtain one of these, take the combinatorial boundary $\Omega$ of $G$, consisting of all reduced one-way
infinite strings in the $u_i$'s and their inverses, topologize it compactly with cylinder sets
$\Omega(s) \ ( =$ set of strings in $\Omega$ beginning with the reduced word $s$) , and put on $\Omega$
a Borel probability measure that is quasi-invariant under the natural left action of $G$. Let $P_1,
\ldots, P_n$ be measurable functions such that 
$$|P_i(\omega)|^2 = \frac{d\mu(u_i^{-1} \omega)}{d \mu(\omega)}.$$
For each $i$, we get a unitary operator $U_i$ on $L^2(\Omega, \mu)$ by setting 
$$(U_i f)(\omega) = P_i(\omega) f(u_i^{-1} \omega),$$ 
and hence a unitary representation $\pi$ --- called a boundary representation --- of $G$ on $L^2(\Omega, \mu)$.  All
boundary representations of $G$ are weakly contained in the left regular representation \cite{Spielberg}, so the states
obtained by composing them with vector states are all reduced. What is needed for our purposes would seem to be
readily available: (1) a probability measure $\mu$ (which can be obtained by prescribing the function $s \mapsto
\mu(\Omega(s))$ from $G$ to $[0,1]$ subject to the obvious compatibility requirement), with the Radon-Nikodym
derivatives $d \mu \circ u_i^{-1}/ d \mu$ in good shape; and (2) measurable functions $q_i$ of modulus one such
that 
$$\sum_i c_i q_i(\omega) \sqrt{\frac{d\mu(u_i^{-1} \omega)}{d \mu(\omega)}} = \lambda$$ 
for $\mu$-almost all $\omega$. If $P_i = q_i \sqrt{d(\mu \circ u_i^{-1})/d\mu}$, and $\pi$  is the resulting
unitary representation of $G$, then the unit constant function ${\sf 1}$ is a $\lambda$-eigenvector for $\sum c_i
\pi(u_i)$, and $\phi = \ < \pi( \cdot) {\sf 1}, {\sf 1}>$ is a state of the type we are seeking.

A diligent search, however, finds only one such choice of $\mu$ and associated $q_i$'s (which turn out to
be $\pm 1$-valued). This apparatus is described in Section \ref{bdryrepns} below. We show there
that $\phi = \ < \pi( \cdot) {\sf 1}, {\sf 1}>$ is of the form $\phi_{B, \vec{c}}$, where the
entries of $B$ are obtained as follows. Let ${\mathbb R}^n_+$ denote the positive orthant of ${\mathbb R}^n$, and let
$D_n$ be the set of points $(s_1, \ldots, s_n)$ in ${\mathbb R}^n_+$ satisfying
$$2 \max s_i - \sum s_i < 1 < \sum s_i .$$
It is shown in Section \ref{mapS} below that the map $S : {\mathbb R}^n_+ \rightarrow D_n$ defined by
$$S(x_1, \ldots , x_n) = \frac{1}{t} (x_1(1 + y_1), \ldots, x_n(1 + y_n)),$$
where
$$t =  \sum_i x_i\ \ \mbox{and for each} \ j , \ y_j = \sum_{i \not= j} x_i \ ,$$ is bijective.
Notice that $(c_1^2/\lambda^2, \ldots, c_n^2/\lambda^2) \in D_n$.  Let
$$(x_1, \ldots, x_n) = S^{-1} (c_1^2/\lambda^2, \ldots, c_n^2/\lambda^2).$$
The matrix $B = \{ \phi(u_i^{-1} u_j) \}$ has entry
$$b_{ij} = - \sqrt{\frac{x_i x_j}{(1+y_i) (1 + y_j)}}$$
for $i \not= j.$ The associated $a_i$'s are given by 
$$a_i = \sqrt{\frac{x_i}{t(1 + y_i)}}.$$
With this particular choice of inputs, the state $\phi_{B, \vec{c}}$ becomes the one described in Section \ref{summary}
above.

\

\section{Pure eigenstates}\label{pure}

Fix positive numbers $c_1, \ldots, c_n$ and a positive (semidefinite) $n \times n$ matrix $B = (b_{ij})$ with
real entries such that $b_{ii} = 1$ for each $i$. With $\vec{c} = (c_1, \ldots, c_n)$, let $\lambda = \sqrt{B
\vec{c} \cdot \vec{c}}$. We assume henceforth that $\lambda > 0.$ Let $\vec{a} = \lambda^{-1} B \vec{c}.$
Define $\phi : G \rightarrow {\mathbb R}$ by
$$\phi(s) = \prod_i a_i^{|s|_i - \sum_j \gamma_{ij}(s)} \prod_{i > j} b_{ij}^{\gamma_{ij}(s)}$$ where $|\cdot |_i$ (resp.
$\gamma_{ij}$) counts the number of occurrences of $u_i$ or $u_i^{-1}$ (resp. $u_i^{-1} u_j$ or $u_j^{-1} u_i$) in a
reduced word in $G$. We will show in this section that $\phi$ is a pure
$\lambda$-eigenstate of $G$ for $\sum_i c_i u_i$. 

The algebraic properties of $\phi$ are easily established.

\

\begin{lem}\label{algprop} (a) The function $\phi$ satisfies: $\phi(1) = 1 $; 

$\phi(s) = \phi(s^{-1})$ for all $s$ in $G$; 

$\phi(u_i s) = a_i \phi(s)$ for $s$ in $G$ not beginning with $u_i^{-1}$ ; 

$\phi(s u_i) = a_i \phi(s)$ for $s$ in $G$ not ending in the inverse of a generator; 

$\phi(s u_j^{-1} u_i) = b_{ij} \phi(s)$ for $s$ not ending in $u_j$.

\noindent(b) The properties in part (a) characterize $\phi$ among real functions on $G$.

\noindent(c) For all $s$ in $G$, we have
$$\sum_i c_i \phi(s u_i) = \lambda \phi(s).$$
\end{lem}

\

\begin{proof} (a) All of the ingredients of the definition are unchanged when $s$ is replaced by
$s^{-1}$. The $\gamma$'s and all of the $|\cdot|$'s except for $|\cdot|_i$ are unchanged when
$s$ not beginning with $u_i^{-1}$ is multiplied on the left by $u_i$. The same is true for right multiplication
by $u_i$ provided $s$ does not end in the inverse of a generator. If $i \not= j$ and $s$ does not end in $u_j$,
then right multiplication of $s$ by $u_j^{-1} u_i$ increases $\gamma_{ij}, \gamma_{ji}, |\cdot|_j,$ and
$|\cdot|_i$ by 1, and leaves all other $\gamma$'s and $|\cdot|$'s unchanged. Thus the exponent of $b_{ij} =
b_{ji}$ in the formula for $\phi$ increases by 1, and all other exponents in the formula are unchanged.

(b) The properties in (a) imply that $\phi(u_i^{\pm 1}) = a_i,$ and furthermore permit calculation of $\phi$
on words of a given length from its values on words of smaller length.

(c) If $s$ does not end in the inverse of a generator, then
$$\sum_i c_i \phi(s u_i) = \left(\sum_i a_i c_i\right) \phi(s) = \frac{B \vec{c} \cdot \vec{c}}{\lambda}
\phi(s) = \lambda \phi(s) .$$
Otherwise, $s = t u_j^{-1}$ for some $t$ not ending in $u_j$ and we have
$$\sum_{i \not= j} c_i \phi(s u_i) + c_j \phi(s u_j) = \left(\sum_{i \not= j} b_{ji} c_i + c_j  \right)
\phi(t) = \lambda a_j  \phi(t) = \lambda \phi(s),$$
where we have used $b_{jj} = 1, B \vec{c} = \lambda \vec{a},$ and $\phi(s) = \phi(s^{-1}) = \phi(u_j t^{-1}) = a_j
\phi(t)$. \end{proof}

\

To begin the construction of the representation that has $\phi$ as a matrix entry, let $G^+$ be the unital
semigroup in $G$ generated by $u_1, \ldots, u_n$, and let $G^+_k$ be the set of group elements in $G^+$ of
length $k$. For $k = 1, 2, \ldots,$ let $A_k$ be the $n^k \times n^k$ matrix with entries indexed by $G^+_k
\times G^+_k$ whose $(s,t)$-entry is $\phi(s^{-1}t).$

\

\begin{lem}\label{mxAk} The matrix $A_k$ is positive for $k = 1, 2, \ldots \ .$ 
\end{lem}

\

\begin{proof} Notice that $A_1 = B.$ For the inductive step, regard $G^+_{k+1}$ as the disjoint union of $n$ 
copies of $G^+_k$ by writing $$G^+_{k+1} = u_1 G^+_k \cup u_2 G^+_k \cup \dots \cup u_n G^+_k \ .$$ We can then
write $A_{k+1}$ in terms of $A_k$ as an $n \times n$ matrix of $n^k \times n^k$ matrices. For $s, t$ in $G^+_k,$
we have 
$$\phi(s^{-1} u_i^{-1} u_j t) = \left\{ \begin{array}{ll} \phi(s^{-1}t) & i = j \\
b_{ij} \phi(s) \phi(t) & i \not= j \end{array} \right. \ .$$
Thus the $i,j$ entry of $A_{k+1}$, viewed in this way, is $A_k$ if $i = j$ and $b_{ij} X_k$ if $i \not= j,$
where  $X_p$, for a positive integer $p$,  is the matrix of the (positive, one-dimensional) operator on
$\ell^2(G^+_p)$ defined by
$$(X_p \vec{\xi})(s) = \left(\sum_{t \in G^+_p} \xi(t) \phi(t) \right) \phi(s).$$
Notice that $X_1 \vec{\xi} = (\vec{\xi} \cdot \vec{a}) \vec{a}$ for $\vec{\xi}$ in $\ell^2(G^+_1).$
We may write
$$A_{k+1} = A_k \otimes I + X_k \otimes (B - I) = (A_k - X_k) \otimes I + X_k \otimes B.$$
in tensor products of $G^+_k \times G^+_k$ matrices and $n \times n$ matrices. Now for $s$ in $G^+$ , we have
$\phi(u_i s) = a_i \phi(s),$ whence it follows that
$X_{k+1} = X_k \otimes X_1.$ This means that
$$A_{k+1} - X_{k+1} = (A_k - X_k) \otimes I + X_k \otimes (B - X_1).$$
Once we show that $B - X_1$ is positive, it will follow by induction that  $A_k - X_k,$ and hence
$A_k$, is positive for all $k$. For $\vec{\xi} = (\xi_1, \ldots, \xi_n)$ in $\ell^2(G^+_1)$, we have
$$X_1 \vec{\xi} \cdot \vec{\xi} = |\vec{\xi} \cdot \vec{a}|^2 = \frac{|\vec{\xi} \cdot B
\vec{c}|^2}{\lambda^2} = \frac{|B^{1/2} \vec{\xi} \cdot B^{1/2} \vec{c}|^2}{\lambda^2} \leq \frac{(B \vec{\xi}
\cdot \vec{\xi}) (B \vec{c} \cdot \vec{c})} {\lambda^2} = B \vec{\xi} \cdot \vec{\xi} , $$
showing $B - X_1 \geq 0$ as promised.  \end{proof}

\

The Hilbert space $H$ of the representation we seek is constructed as follows. By Lemma
\ref{mxAk}, there is for each positive integer $k$  a finite dimensional complex inner
product space  $E_k$ spanned by vectors $\{\Delta_s : s \in G^+_k \}$ with inner
product $< \cdot, \cdot >$ satisfying $<\Delta_t , \Delta_s> = \phi(s^{-1}t).$ (We write 
$E_0$ for the one-dimensional inner product space spanned by the unit vector
$\Delta_1.$) Because
$$ \lambda^{-2} \sum_{i,j} c_i c_j \phi(u_j^{-1} t^{-1} s u_i) = \phi(t^{-1}s)$$
for all $s, t$ in $G^+$, we have an isometry from $E_k$ into $E_{k+1}$
for each $k$ sending $\Delta_s$, for $s$ in $G^+_k$,  to $\lambda^{-1} \sum_i c_i \Delta_{su_i}.$ Let $H_0$ be
the Hilbert space inductive limit of the resulting tower $E_0 \rightarrow E_1  \rightarrow E_2
\rightarrow \ldots .$ Then $H_0$ is the closed linear span of $\{\Delta_s : s \in
G^+ \}$, and these vectors satisfy 
$$\sum_i c_i \Delta_{su_i} = \lambda \Delta_s \ \ \mbox{and} \ \  <\Delta_t , \Delta_s> \ = \
\phi(s^{-1}t).$$  
Left multiplication by each generator $u_i$ gives rise to an isometry $V_i$ of 
$H_0$ into itself. Let $H_i' = H_0 \ominus V_i H_0.$  For each $i,$ let
$S^-_i$ be the subset of $G$ consisting of the reduced words ending in 
$u_i^{-1}$. The Hilbert space
$H$ is 
$$H = H_0 \oplus \bigoplus_{i = 1}^n \left(\ell^2 (S^-_i) \otimes H_i'
\right) .$$  
For each $i$, let $U_i$ be the unitary operator on $H$ that maps $H_0$ to $V_i H_0$ by $V_i$, maps
$\delta_{u_i^{-1}} \otimes H_i'$ to $H_i' =  H_0 \ominus V_iH_0$ by erasing the tensor, and maps
$\delta_s
\otimes \eta$ to $\delta_{u_i s} \otimes \eta$ for all other $s$ ending in the inverse of a
generator, and for $\eta$ in the appropriate space $H_j'$. Denote by $\pi$ the
unitary representation of $G$ on $H$ that takes $u_i$ to $U_i$. Our next goal is to show that
$\phi = \ <\pi(\cdot) \Delta_1, \Delta_1>$.

\ 

Write $P_i$ for the orthogonal projection of $H_0$ on $H_i'$.

\

\begin{lem}\label{orthdec}  \ \ \ (a) $P_i \Delta_1 = \Delta_1 - a_i \Delta_{u_i}$ and $U_i^* \Delta_1 = \delta_{u_i^{-1}}
\otimes P_i \Delta_1 + a_i \Delta_1$.

\noindent (b) For $s$ in $u_jG^+$ and $i \not= j,$ we have $P_i \Delta_s = \Delta_s - \phi(u_i^{-1} s) \Delta_{u_i}$ and
$$U_i^* \Delta_s = \delta_{u_i^{-1}} \otimes P_i \Delta_s + \phi(u_i^{-1} s) \Delta_1 .$$

\noindent (c) For $\eta \in U_j H_0$ and $i \not= j,$ we have $(1 - P_i) \eta = \ <\eta, \Delta_{u_i}> \Delta_{u_i}$.
\end{lem}

\

\begin{proof} (a) Both statements follow from the observation that
$$<\Delta_{u_i t}, \Delta_1 - a_i \Delta_{u_i}>  \ = \phi(u_i t) - a_i \phi(t) = 0$$
for all $t$ in $G^+$.

(b) Notice that
$$<\Delta_s - \phi(u_i^{-1} s) \Delta_{u_i}, \Delta_{u_i t}> \ = \phi(t^{-1} u_i^{-1} s) - \phi(t^{-1}) \phi(u_i^{-1} s) = 0$$
for all $t$ in $G^+$ because there is no cancellation in $u_i^{-1} s$

(c) The $\Delta_s$'s with $s$ as in (b) span a dense subspace of $U_j H_0$, so this follows from (b). \end{proof}

\

\begin{lem}\label{mapdel} For each $t$ in $G$, the vector $\pi(t) \Delta_1$ may be written in the form $\xi + h
\Delta_w,$ where $h$ is a real number, $w \in G^+,$ and $\xi$ is either zero or a finite sum of terms $\delta_r
\otimes
\eta_r,$ where each $r$ belongs to some $S_i^-$ and begins with the same generator or inverse generator as $t$,
and $\eta_r$ belongs to the corresponding $H_i'$. Furthermore, $w = 1$ if $t$ begins with the inverse of a
generator.
\end{lem}

\

\begin{proof} This is obvious if $t \in G^+$ (the case in which $\xi = 0).$ By what we have observed
just above, $t = u_i^{-1}$ makes $\xi = \delta_{u_i^{-1}} \otimes P_i \Delta_i, \ h = a_i,$ and $w = 1.$ For $s$
in $u_j G^+$ and $i \not =j,$ we also have
$$\pi(u_i^{-1} s) \Delta_1 = \delta_{u_i^{-1}} \otimes P_i \Delta_s + \phi(u_i^{-1} s) \Delta_1,$$
and furthermore 
$$\pi(u_k^{-1} u_i^{-1} s) \Delta_1 = \delta_{u_k^{-1} u_i^{-1}} \otimes P_i \Delta_s + \phi(u_i^{-1} s)
(\delta_{u_k^{-1}} \otimes P_k \Delta_1 + a_k \Delta_1),$$
while 
$$\pi(u_k u_i^{-1} s) \Delta_1 = \delta_{u_k u_i^{-1}} \otimes P_i \Delta_s + \phi(u_i^{-1} s)
\Delta_{u_k}$$
when $k \not= i.$ And so forth --- the asserted form plainly remains intact under further noncancelling left
multiplication by generators or inverse generators. \end{proof}

\

\begin{prop}\label{vecstate} For all $t$ in $G$, we have $<\pi(t) \Delta_1, \Delta_1> \ = \ \phi(t).$
\end{prop}

\

\begin{proof} Write $\pi(t) \Delta_1 = \xi + h \Delta_w$ as in Lemma \ref{mapdel}.  If $t$ does not begin
with $u_i^{-1},$ then $U_i \xi \perp H_0$ and so
$$<\pi(u_i t) \Delta_1, \Delta_1> \ = \ h <\Delta_{u_i w}, \Delta_1> \ =  h a_i \phi(w) = \ a_i < \pi(t) \Delta_1, \Delta_1>
.$$ Furthermore, for $j \not= i$ we have $U_j^* U_i \xi \perp H_0$ and
$$<\pi(u_j^{-1} u_i t) \Delta_1, \Delta_1> \ = \ h < \Delta_{u_i w}, \Delta_{u_j}>$$ $$= h \phi(u_j^{-1} u_i
w) = b_{ij} h \phi(w) = b_{ij} < \pi(t) \Delta_1, \Delta_1> .$$
If $t$ begins with $u_i^{-1}$, then $w = 1$, and for any $j$
$$<\pi(u_j^{-1} t) \Delta_1, \Delta_1> \ = \ <\pi(t) \Delta_1, \Delta_{u_j}> \ = h <\Delta_1, \Delta_{u_j}>$$
$$= a_j h = a_j < \pi(t) \Delta_1, \Delta_1>.$$ 
Notice also that $<\pi(t) \Delta_1, \Delta_1> \ = \ <\pi(t^{-1}) \Delta_1, \Delta_1>$  (because $h$ is real). It now follows
from part (b)  of Lemma \ref{algprop}  that $<\pi(\cdot) \Delta_1, \Delta_1>$ coincides with $\phi.$ \end{proof}

\

\begin{thm}\label{irred} The representation $\pi$ is irreducible.
\end{thm}

\

\begin{proof} We begin by noticing that $\Delta_1$ is cyclic  for $\pi$. Indeed, by
Lemma \ref{orthdec}, the linear span of $\pi(G) \Delta_1$ contains $\delta_{u_i^{-1}}
\otimes P_i \Delta_s$ for each $s$ in $G^+$ and each $i$. It follows that the closed linear span contains each
summand $\ell^2(S_i^-) \otimes H_i'$. Of course, it also contains $H_0$ because $\pi(s) \Delta_1 = \Delta_s$ for $s$
in $G^+$.

We will show in several steps that the commutant $\pi(G)'$ of $\pi(G)$ consists of scalars. Our argument makes
essential use of a result of P. Linnell (Proposition 1.4 and especially Lemma 3.6 of \cite{Linnell}) which implies that a
nonzero linear combination of elements of $G$ must have kernel zero in the left regular representation.

\underline{Claim 1} \ If $T \in \pi(G)'$, then $T \Delta_1$ is orthogonal to $\ell^2(S_i^-) \otimes P_i
\Delta_1$ for each $i$.

\underline{Proof of Claim 1} \ Since $\Delta_1$ is a $\lambda$-eigenvector for $\sum_i c_i \pi(u_i),$ so is $T
\Delta_1.$ Thus 
$$\lambda <T \Delta_1, \Delta_1> \ = \ \sum_i c_i<T \Delta_1, U_i^* \Delta_1>$$
$$= \sum_i c_i <T \Delta_1, \delta_{u_i^{-1}} \otimes P_i \Delta_1 + a_i \Delta_1>$$
$$= \ \lambda <T \Delta_1, \Delta_1> + \sum_i c_i <T \Delta_1, \delta_{u_i^{-1}} \otimes P_i
\Delta_1> ,$$
where we have used Lemma \ref{orthdec} and the identity $\sum_i a_i c_i = \lambda$.
This makes
$$(*) \ \ \sum_i c_i <T \Delta_1, \delta_{u_i^{-1}} \otimes P_i \Delta_1> \ = 0 .$$
Define $\xi$ on $G$ by setting $\xi(t) = \ <T \Delta_1, \ \delta_t \otimes P_i \Delta_1>$ if $t \in S_i^-$
for some $i$, and $\xi(t) = 0$ for all other $t$ in $G$. Then
$$\sum_j c_j \xi(u_j^{-1} t) = \lambda \xi(t)$$
(1) for $t$ in  $\cup_i S_i^-$ because $T \Delta_1$ is a $\lambda$-eigenvector for $\sum_j c_j \pi(u_j)$;
(2) for $t = 1$ by (*) above, because $\xi(1) = 0$; and (3) for all other $t$ in $G$ because for such $t$ all
of the values of $\xi$ in the asserted relationship are zero. Because the vectors $\delta_t \otimes P_i
\Delta_1$ that appear in the nonzero part of the definition of $\xi$ are mutually orthogonal and norm-bounded,
the function $\xi$ belongs to $\ell^2(G).$ The claim now follows from 1.4 in \cite{Linnell}.

\underline{Claim 2} \ For $T$ in $\pi(G)'$ and $s$ in $G^+,$ we have 
$$<T \Delta_1, \Delta_{u_i s}> \ = \ a_i <T \Delta_1, \Delta_s>$$ 
for each $i$, and hence
$$<T \Delta_1, \Delta_s> \ = \ \phi(s) <T \Delta_1, \Delta_1>$$
for all $s$ in $G^+$.

\underline{Proof of Claim 2} \ If furthermore $T = T^*$, we see using Lemma \ref{orthdec} and Claim 1 that
$$<T \Delta_1, \Delta_{u_i s}> \ = \ <\delta_{u_i^{-1}} \otimes P_i \Delta_1 + a_i \Delta_1, \pi(s) T
\Delta_1>$$
$$ = \ <\delta_{s^{-1} u_i^{-1}} \otimes P_i \Delta_1, T \Delta_1> + a_i <\Delta_1, T \Delta_s> \ = \ a_i <T
\Delta_1, \Delta_s> \ .$$
Because the claimed relationship is linear in $T$, we can remove the assumption that $T$ be selfadjoint. The
second assertion follows by induction on the length of $s$.

\underline{Claim 3} \ $T \Delta_1 \in H_0$ for $T$ in $\pi(G)'$.

\underline{Proof of Claim 3} \ It will suffice to show that 
$$T \Delta_1 \perp \ell^2(S_j^-) \otimes P_j \Delta_{u_i s}$$
for $j \not= i$ and $s$ in $G^+$. We proceed more or less as in the proof of Claim 1. Fix $s$ in $G^+$ and the
index $i$. Define $\xi$ on $G$ by setting $\xi(t) = \ <T \Delta_1, \delta_t \otimes P_j \Delta_{u_i s}>$ if $t
\in S_j^-$ for some $j \not= i,$ and $\xi(t) = 0$ for all other $t$ in $G$. Thus $\xi \in \ell^2(G).$ As in the proof of
Claim 1, we will show that $\xi$ must vanish by showing that it satisfies
$$\sum_k c_k \xi(u_k^{-1}t) = \lambda \xi(t)$$
for all $t$ in $G$. This identity holds if $t \in S_i^-$ or if $t$ ends in a positive generator power because
all of the values of $\xi$ that appear are zero. It holds if $t \in S_j^-$ for some $j \not= i$ because $T
\Delta_1$ is a $\lambda$-eigenvector for $\sum_k c_k \pi(u_k)$. The only remaining case is $t = 1,$ for which
we must show that
$$\sum_{k \not= i} c_k \xi(u_k^{-1}) = 0 \ .$$
But
$$\lambda a_i \phi(s) <T \Delta_1, \Delta_1> \ = \lambda <T \Delta_1, \Delta_{u_i s}> \ =
\sum_k c_k<T \Delta_1, U_k^* \Delta_{u_i s}>$$
$$= c_i <T \Delta_1, \Delta_s> + \sum_{k \not= i} c_k <T \Delta_1, \delta_{u_k^{-1}} \otimes P_k \Delta_{u_i s}
+ b_{ik} \phi(s) \Delta_1>$$
$$= c_i \phi(s) <T \Delta_1, \Delta_1> + \sum_{k \not= i} b_{ik} c_k \phi(s) <T \Delta_1, \Delta_1> \ + \sum_{k
\not= i} c_k \xi(u_k^{-1}) $$
$$= \lambda a_i \phi(s) <T \Delta_1, \Delta_1> + \sum_{k \not= i} c_k \xi(u_k^{-1}),$$
where we have used Claim 2, Lemma \ref{orthdec} above, and $B \vec{c} =\lambda
\vec{a}$. It follows from Claim 2 that $T \Delta_1 - <T \Delta_1, \Delta_1> \Delta_1$ is orthogonal to $\Delta_s$
for every $s$ in $G^+$, and hence $T \Delta_1 = \ <T \Delta_1, \Delta_1> \Delta_1$ by Claim 3, for all $T$ in
$\pi(G)'$. Because $\Delta_1$ is cyclic for $\pi$, this proves that $\pi$ is irreducible. \end{proof}

\

\section{Boundary representations}\label{bdryrepns}

\

The representations considered in this section are all weakly contained in the left regular representation,
so we must restrict attention to $\lambda$'s in the spectrum of $\sum c_i u_i$ in that representation.

\

\begin{prop}\label{redspectrum} Let complex coefficients $c_1, c_2, \ldots, c_n$ be given, and let $X = \sum c_i u_i$
The reduced spectrum of $X$  is 
$$\{\lambda \in {\mathbb C} : r_0 \leq |\lambda| \leq \left(\sum |c_i|^2 \right)^{1/2}\},$$
where $r_0^2$ is zero if none of the $|c_i|^2$'s exceeds the sum of the others, and otherwise the maximum of the
$|c_i|^2$'s minus the sum of the others.
\end{prop}

\

\begin{proof} The spectrum is connected because the reduced $C^*$-algebra of $G$ contains no nontrivial
idempotents \cite{PimsnerVoiculescu}. It is rotationally invariant because there is an automorphism of this
$C^*$-algebra that multiplies each $u_i$ by a given scalar of modulus one. Thus, the spectrum must be either a closed
disc about 0 or a closed annulus. Let $|| \cdot ||_{op}$ denote the  reduced operator norm. Then for every positive
integer $k$ we have
$$||X^k||_2 \leq ||X^k||_{op} \leq (k+1)||X^k||_2,$$
where the lower bound is obvious and the upper bound follows from Haagerup's inequality \cite{Haagerup} plus the
observation that $X^k$ is a linear combination of words of length $k$.  Since the coefficient of a $k$-fold product of
generators in $X^k$ is the corresponding product of $c$'s, we have
$$||X^k||_2 = ||X||_2^k =  \left(\sum |c_i|^2 \right)^{k/2},$$
so
$$\left(\sum |c_i|^2 \right)^{1/2} \leq ||X^k||_{op}^{1/k} \leq (k+1)^{1/k} \left(\sum |c_i|^2 \right)^{1/2}.$$
This shows that the reduced spectral radius of $X$ is the 2-norm of the coefficient vector.  It remains to show that we
have correctly identified the inner radius $r_0,$ which we do by induction on the number of nonzero coefficients. In
the case of only one nonzero coefficient, say $c_1$, the spectrum is the circle about 0 of radius $|c_1|,$ which plainly
coincides with $r_0$. Suppose that the proposition gives the correct inner radius when there are  $m-1$
nonzero coefficients. Without loss of generality, we may assume $X = c_1 u_1 + \ldots + c_{m-1} u_{m-1} + u_m,$
where $|c_1|, \ldots, |c_{m-1}| \leq 1.$ Let $Y = (c_1 u_1 + \ldots c_{m-1} u_{m-1}) u_m^{-1}$. There are two
cases: (1) $|c_1|^2 + \ldots + |c_{m-1}|^2 \geq 1$; and (2) the contrary. In case (1), the inner radius of the reduced
spectrum of $Y$ is at most 1 by our induction hypothesis (plus the assumption that the $|c|$'s are all at most 1) and
the outer radius is at least 1 by the spectral radius formula already established. Thus $X  \ (= (Y + 1) u_m)$ is not
invertible in the left regular representation, giving $r_0 = 0$ as required in this case.  In case (2), let
$$r = \sqrt{1 - |c_1|^2 - \ldots -|c_{m-1}|^2}.$$
If $|\lambda| < r,$ then the reduced spectral radius of $Y - \lambda u_m^{-1}$ is less than 1,
and hence $ X - \lambda \ (= (Y - \lambda u_m^{-1} + 1)u_m)$ is invertible in the left regular representation. This
shows that the reduced spectrum of $X$ has inner radius at least $r$. On the other hand, $X - r$ is singular in the left
regular representation because $Y - r u_m^{-1}$ has reduced spectral radius precisely 1. This shows that the reduced
spectrum of $X$ has the correct inner radius in case (2), completing the proof. \end{proof}

\

Here is a review of some essential facts about the combinatorial boundary of $G$.  Let $V$ be the set $\{u_1, u_1^{-1},
u_2, u_2^{-1},\ldots, u_n, u_n^{-1}\}.$ The combinatorial boundary
$\Omega$ of $G$ is the set of all infinite strings $\omega = v_1 v_2 \ldots $ where each $v_i \in V$ and $v_{i+1}
\not= v_i^{-1}$ for all $i$. It becomes a compact Hausdorff space when equipped with the product topology, a basis
for which consists of cylinder sets described as follows. For each positive integer $k$, let $p_k$ be the map
from $\Omega$ to $G$ that reads the first $k$ symbols. For $s$ in $G$ with $|s| = k,$ let $\Omega(s) =
p_k^{-1}(\{s\}),$ that is, the set of strings in $\Omega$ beginning with $s$. Also let $\Omega(1) = \Omega.$ The
$\Omega(s)$'s are all open and closed, and form a basis for the topology of $\Omega.$ For each $s$ in $G$, let
$G(s)$ be the set of reduced words in $G$ beginning with $s$, with $E_1 = G$.  Further let $E_s$ be the orthogonal
projection of $\ell^2(G)$ on $\ell^2(G(s))$, and let $A$ be the $C^*$-algebra generated by the $E_s$'s --- that is, the
norm closure of their linear span. Then $A$ is commutative, and $A$ modulo its intersection with the ideal $K$ of
compact operators is easily seen to be isomorphic to $C(\Omega)$. (Let $\phi$ be a
multiplicative linear functional on $A$ that kills $A \cap K.$ Since each $E_s$ is  a one-dimensional operator
plus the sum of $E_{sv}$ over $v$ in $V$ such that $sv$ is reduced, there is a unique $\omega$ in
$\Omega$ such that $\phi(E_{p_m(\omega)}) = 1$ for all $m$. On the other hand, given $\omega$ in $\Omega$, one
obtains a multiplicative linear functional on $A$ by taking a weak$^*$ limit of the vector states $<\cdot \ \
\delta_{p_k(\omega)}, \delta_{p_k(\omega)}>$. In this way, one has a continuous bijection between $\Omega$ and
the maximal ideal space of $A/(A \cap K)$.) Borel measures on $\Omega$ are determined by their values on cylinder
sets, and can be defined by prescribing those values subject only to an obvious compatibility condition.  In fact, given
any function
$\hat{\mu} : G \rightarrow [0,1]$ with $\hat{\mu}(1) = 1$ such that
$$\hat{\mu}(s) = \sum \{\hat{\mu}(sv) : v \in V, |sv| = |s| + 1\}$$
for every $s$ in $G$, there is a Borel probability measure $\mu$ such that $\mu(\Omega(s)) = \hat{\mu}(s)$ for
all $s$. (To see this, notice that for every $k$, the sum of the values of $\hat{\mu}$ on the words of length $k$ is
1. Let $\xi_k$ be $\sqrt{\hat{\mu}}$ times the indicator function of the set of words of length $k$, so $\xi_k$ is a
unit vector in  $\ell^2(G)$. Any weak$^*$-limit of the corresponding vector states on $A$ gives the desired $\mu$,
because $< E_s \xi_k, \xi_k > = \hat{\mu}(s)$ for $k \geq |s|.$) We will continue to use the notation $\hat{\mu}$
for the function on $G$ corresponding to  measure $\mu$ on $\Omega.$ 

One type of eligible $\hat{\mu}$ can be specified by choosing functions $\beta: V \rightarrow [0,1]$ and
$\alpha : V \times V \rightarrow [0,1]$ such that
$$\sum_{v \in V} \beta(v) = 1 \ , \sum_{v' \in V} \alpha(v,v') = 1 \ \ \forall v \in V \ , \ \alpha(v,v^{-1}) = 0
\ \ \forall v \in V \ ,$$
and then setting $\hat{\mu}(1) = 1, \ \hat{\mu}(v) = \beta(v)$ for $v$ in $V$, and for reduced words of length two
or greater, \ $\hat{\mu}(v_1 v_2 \ldots v_k) = \beta(v_1) \alpha(v_1, v_2) \ldots \alpha(v_{k-1},v_k).$ (In other
words, the measure $\mu$ is  the probability measure on the space of sample paths in a Markov chain with states
labeled by $V$, transition probabilities given by $\alpha,$ and initial probabilities given by $\beta$.) If all of the
values of $\beta$ and of $\alpha$ (except the $\alpha(v,v^{-1})$'s) are positive, the measure
$\mu$ is quasi-invariant under the natural left action of $G$ on $\Omega$, and the Radon-Nikodym derivatives of the
translates of $\mu$ by the generators of $G$ are easily calculated. Namely we have
$$\frac{d\mu(u_i^{-1} \omega)}{d\mu(\omega)} = \left\{\begin{array}{ll} \frac{\beta(v)}{\beta(u_i)
\alpha(u_i,v)}
\ \ & \ \ p_2(\omega) = u_i v \\&\\ \frac{\beta(u_i^{-1}) \alpha(u_i^{-1},v)}{\beta(v)} \ \ & \ \ p_1(\omega) =
v \not= u_i \end{array} \right. \ $$
because for reduced words $s$ of length at least 2, the ratio $\mu(u_i^{-1} \Omega(s))/\mu(\Omega(s))$ depends
only on the first two symbols in $s$ (in the manner indicated by the switches in the formula).

Let positive numbers $c_1, c_2, \ldots, c_n$ be given, as well as a positive number $\lambda$ satisfying
$$\lambda^2 < \sum_i c_i^2 \ \ \mbox{and} \ \ c_j^2 - \sum_{i \not= j} c_i^2  < \lambda^2  \ \ \mbox{for} \ \ j =
1, \ldots, n \ .$$ Our immediate aim is to exhibit a Borel probability measure $\mu$ on $\Omega$ and a unitary
representation of $G$ on $L^2(\Omega,\mu)$ in which $c_1 u_1 + c_2 u_2 + \ldots c_n u_n$ has $\lambda$ as an
eigenvalue. This can be done fairly cleanly in terms of the inverse of the map $S$ defined on the positive orthant
of ${\mathbb R}^n$ by
$$S(x_1, x_2, \ldots, x_n) = \left( \frac{x_j (1 + \sum_{i \not= j} x_i)}{\sum_i x_i} \right)_{j=1}^n .$$
It is shown in the appendix that $S$ is injective on the positive orthant. Furthermore, the conditions we have
imposed on $\lambda$ ensure that $$\lambda^{-2} (c_1^2, c_2^2, \ldots, c_n^2)$$ belongs to the range of $S$.
Accordingly, we write
$$(x_1, x_2, \ldots, x_n) = S^{-1} \left( \frac{c_1^2}{\lambda^2}, \frac{c_2^2}{\lambda^2}, \ldots,
\frac{c_n^2}{\lambda^2} \right) \ .$$
As in the appendix, we write $t = \sum_i x_i$ and $y_j = t - x_j = \sum_{i \not = j} x_i.$ The measure $\mu$ is
the one constructed as in the previous paragraph for $\beta$ and $\alpha$ defined by
$$\beta(u_i) = \frac{x_i}{t (1 + t)} \ , \ \beta(u_i^{-1}) = \frac{x_i}{1 + t} ,$$
$$\alpha(u_j, u_i) = \alpha(u_j^{-1}, u_i^{-1}) = \frac{x_i}{t (1 + y_j)} ,$$
$$\alpha(u_j, u_i^{-1}) = \alpha(u_j^{-1}, u_i) = \frac{x_i}{1 + y_j} \ \ (i \not= j) ,$$
and $\alpha(u_i,u_i^{-1}) = 0  = \alpha(u_i^{-1},u_i).$ It is readily checked that $\alpha$ and $\beta$ satisfy
all of the sum-to-1 conditions of the previous paragraph. Radon-Nikodym derivatives under translation by the
generators are given by
$$\frac{d\mu(u_i^{-1} \omega)}{d\mu(\omega)} = \left\{\begin{array}{ll} \frac{t(1 + y_i)}{x_i} \ \ & \ \
p_1(\omega) = u_i \\&\\ \frac{t x_i}{1 + y_i} \ \ & \ \ p_1(\omega) = u_j , \ j \not= i \\&\\ 
\frac{x_i}{t(1+y_i)} \ \ & \ \ p_1(\omega) = u_j^{-1} \end{array} \right. \ \ .$$
Define $P : \{u_1, u_2, \ldots, u_n\} \times \Omega \rightarrow R$ by
$$P(u_i, \omega) = \left\{\begin{array}{ll} \sqrt{\frac{t(1 + y_i)}{x_i}} \ \ & \ \
p_1(\omega) = u_i \\&\\ -\sqrt{\frac{t x_i}{1 + y_i}} \ \ & \ \ p_1(\omega) = u_j , \ j \not= i \\&\\ 
\sqrt{\frac{x_i}{t(1+y_i)}} \ \ & \ \ p_1(\omega) = u_j^{-1} \end{array} \right. \ \ $$
so $P(u_i, \omega)^2 = d\mu(u_i^{-1} \omega)/d\mu(\omega)$. Define unitaries $U_1, U_2, \ldots, U_n$ on
$L^2(\Omega, \mu)$ by
$$(U_i \xi)(\omega) = P(u_i,\omega) \xi(u_i^{-1} \omega) \ ,$$
and let $\pi$ be the unitary representation of $G$ on this Hilbert space that sends each $u_i$ to the
corresponding $U_i$. It follows from Theorem 2.7 in \cite{Spielberg} (see also Theorem 1X in \cite{KuhnSteger})
that $\pi$ is weakly contained in the left regular representation.

\

\begin{prop}\label{rightsum} We have 
$$\sum_i c_i \pi(u_i) {\sf 1} = \lambda {\sf 1} ,$$
where $\sf 1$ is the unit constant function on $\Omega$.
\end{prop}

\

\begin{proof} Since $\pi(u_i){\sf 1} = P(u_i, \cdot),$ this amounts to showing that
$$\sum_i \frac{c_i}{\lambda} P(u_i,\omega) = 1$$
for all $\omega.$ By construction,
$$\frac{c_i}{\lambda} = \sqrt{\frac{x_i (1 + y_i)}{t}}$$
for each $i$. If $\omega \in \Omega(u_j^{-1})$ for some $j$, we have
$$\sum_i \frac{c_i}{\lambda} P(u_i,\omega) = \sum_i \sqrt{\frac{x_i (1 + y_i)}{t}} \sqrt{\frac{x_i}{t (1 +
y_i)}} = \sum_i \frac{x_i}{t} = 1.$$
Otherwise, $\omega \in \Omega(u_j)$ for some $j$, and we have
$$\sum_i \frac{c_i}{\lambda} P(u_i,\omega) = -\sum_{i \not= j} \sqrt{\frac{x_i (1 + y_i)}{t}} \sqrt{\frac{t x_i}{
1 + y_i}} + \sqrt{\frac{x_j (1 + y_j)}{t}} \sqrt{\frac{t (1 + y_j)}{x_j}}$$
$$ = - \sum_{i \not= j} x_i + 1 + y_j = 1.$$ \end{proof}

\

We now take up the project of identifying the reduced state $<\pi( \cdot) {\sf 1, 1}>.$ Let us call this state
$\phi.$ To begin with, $\phi$ can be expressed in terms of the cocycle $P$ (with values in the multiplicative group of
nonzero reals) on $G \times \Omega$ that extends the function $P$ defined above. Thus 
$$(\pi(s) \xi)(w) = P(s,w) \xi(s^{-1} \omega)$$
for $s$ in $G$ and $\omega$ in $\Omega.$ In particular, $(\pi(s) {\sf 1})(\omega) = P(s,\omega)$ and so
$$\phi(s) = \int_\Omega P(s,\omega) d\mu(\omega) .$$
The cocycle identity satisfied by $P$ is
$$P(r s, \omega) = P(r, \omega) P(s, r^{-1} \omega).$$
Since $P(1,\omega) = 1,$ we have in particular
$$P(u_i^{-1}, \omega) = \frac{1}{P(u_i, u_i \omega)} = 
\left\{\begin{array}{ll} \sqrt{\frac{x_i}{t(1 + y_i)}} \ \ &
\ \ p_1(\omega) \not= u_i^{-1} \\&\\ -\sqrt{\frac{1 + y_i}{t x_i}} \ \ & \ \ p_2(\omega) = u_i^{-1} u_j 
\\&\\ 
\sqrt{\frac{t(1+y_i)}{x_i}} \ \ & \ \ p_2(\omega) = u_i^{-1} u_j^{-1} \end{array} \right. \ \ .$$
It also follows from the cocycle identity that  for an arbitrary reduced word $v_1 v_2 \ldots v_k$ with each $v_i$ in
$V$ we may write $P(v_1 v_2 \ldots v_k, \omega)$ as
$$P(v_1, \omega) P(v_2, v_1^{-1} \omega) P(v_3, v_2^{-1} v_1^{-1} \omega) \ldots P(v_k, v_{k-1}^{-1} \ldots
v_2^{-1} v_1^{-1} \omega).$$ 
This formula leads to the following useful observation.

\

\begin{lem}\label{constant} \ \ \ (a) If the reduced word $r$ begins with $u_i^{-1}$ for some $i$, then $P(r,
\cdot)$ is constant on $\Omega(u_1) \cup \Omega(u_2) \cup \ldots \cup \Omega(u_n).$

(b) On the other hand, $P(r, \cdot)$ is constant on $\Omega(u_i^{-1})$ provided $r$ does not begin with $u_i^{-1}$.
\end{lem}

\

\begin{proof} (a) Write $\Omega(+) = \Omega(u_1) \cup \Omega(u_2) \cup \ldots \Omega(u_n),$ and $r =
u_i^{-1} v_2 \ldots v_k.$ Observe that $P(u_i^{-1}, \cdot)$ is constant on $\Omega(+)$, in fact on $\Omega \setminus
\Omega(u_i^{-1}).$ Whether $v_2$ is $u_j$ for some $j \not= i$ or $u_j^{-1}$, the factor $P(v_2, u_i (\cdot))$ in the
formula is constant on $\Omega(+)$; notice that $u_i \Omega(+) \subseteq \Omega(u_i)$. The remaining factors in the
formula for $P(r, \cdot)$ are constant on $\Omega(+)$ because for $v$ in $V$, the function $P(v, \cdot)$ reads at most
only the first two symbols in its argument.

(b) Suppose $r = u_j^{-1} v_2 \ldots v_k$ for some $j \not= i$ (and hence $v_2 \not= u_j$).  Observe that $P(u_j^{-1},
\cdot)$ is constant on $\Omega(u_i^{-1})$, in fact on $\Omega \setminus \Omega(u_j^{-1})$. We have $u_j
\Omega(u_i^{-1}) \subseteq \Omega(u_j u_i^{-1}),$ so $P(v_2, u_j (\cdot))$ is constant on $\Omega(u_i^{-1})$, and the
remaining factors are constant there for the same reason as in part (a). Likewise if $r = u_j v_2 \ldots v_k$ (and hence $v_2 \not= u_j^{-1}$), the two
initial factors $P(u_j, \cdot)$ and $P(v_2, u_j^{-1} (\cdot))$ are constant on $\Omega(u_i^{-1})$. The remaining
factors are constant on $\Omega(u_i^{-1})$ for the same reason as in part (a). \end{proof}

\

We can now undertake the  calculation that exhibits $\phi$ as one of the states investigated in
Section \ref{pure}

\

\begin{prop}\label{stateformula} For each $i$ and each $j$ different from $i$, let
$$a_i = \sqrt{\frac{x_i}{t(1 + y_i)}} \ \ \mbox{and} \ \ b_{ij} = - \sqrt{\frac{x_i x_j}{(1 + y_i) (1 + y_j)}}
\ .$$ Then 
$$\phi(s) = \prod_i a_i^{|s|_i - \sum_j \gamma_{ij}(s)} \prod_{i > j} b_{ij}^{\gamma_{ij}(s)}$$
where $| \cdot |_i$ (resp. $\gamma_{ij}$) counts the number of occurrences of $u_i$ or $u_i^{-1}$ (resp.
$u_i^{-1} u_j$ or $u_j^{-1} u_i$) in a reduced word in $G$. 
\end{prop}

\

\begin{proof} The calculation is in three parts.

\

1. We show first that if $j \not= i$, and $s$ is a reduced word not beginning with
$u_i^{-1}$ (including the possibility that $s = 1$), then $\phi(u_j^{-1}u_i s) = b_{ij} \phi(s)$.  We have
$$\phi(u_j^{-1} u_i s) = \ < \pi(s) {\sf 1}, \pi(u_i^{-1} u_j) {\sf 1} > \ = \int_{\Omega} P(s,\omega)
\frac{P(u_j,u_i \omega)}{P(u_i,u_i \omega)} d\mu(\omega)$$
$$ = \int_{\Omega \setminus \Omega(u_i^{-1})} + \int_{\Omega(u_i^{-1} u_j)} + \sum_{k \not= i,j}
\int_{\Omega(u_i^{-1} u_k)} + $$
$$+ \sum_k \int_{\Omega(u_i^{-1} u_k^{-1})} P(s,\omega) \frac{P(u_j,u_i
\omega)}{P(u_i,u_i \omega)} d\mu(\omega) .$$
By Lemma \ref{constant}(b), there is a number $C$ such that $P(s,\omega) = C$ for all $\omega$ in
$\Omega(u_i^{-1}).$ Taking this into account, and looking up values for $P(u_i, \cdot )$ and $P(u_j, \cdot )$ in the
various cases that arise, we obtain
$$\phi(u_j^{-1} u_i s) =  - \ \sqrt{\frac{t x_j}{1 + y_j}} \sqrt{\frac{x_i}{t (1 + y_i)}} \int_{\Omega \setminus
\Omega(u_i^{-1})} P(s, \omega) d\mu(\omega)$$
$$ - \ C \sqrt{\frac{t (1 + y_j)}{x_j}} \sqrt{\frac{1 + y_i}{t x_i}}
\hat{\mu}(u_i^{-1} u_j) \ + $$
$$+ \ C \sum_{k \not= i,j} \sqrt{\frac{t x_j}{1 + y_j}} \sqrt{\frac{1 + y_i}{t x_i}} \hat{\mu}(u_i^{-1} u_k) \ +$$
$$ + \ C \sum_k \sqrt{\frac{x_j}{t (1 + y_j)}} \sqrt{\frac{t (1 + y_i)}{x_i}} \hat{\mu}(u_i^{-1} u_k^{-1}) \ .$$
Using $\hat{\mu}(u_i^{-1} u_k^{\pm 1}) = \beta(u_i^{-1}) \alpha(u_i^{-1}, u_k^{\pm 1})$, this makes
$$\phi(u_j^{-1} u_i s) = b_{ij} \int_{\Omega \setminus \Omega(u_i^{-1})} P(s, \omega) d\mu(\omega) \ +$$
$$C \left( -\sqrt{\frac{(1+y_j) (1+y_i)}{x_j x_i}} \cdot \frac{x_i}{1+t} \cdot \frac{x_j}{1 + y_i} + 
\sqrt{\frac{x_j(1 + y_i)}{(1 + y_j) x_i}} \cdot \frac{x_i}{1+t} \sum_{k \not= i,j} \frac{x_k}{1+y_i} + \right.$$
$$ \left. + \sqrt{\frac{x_j (1 + y_i)}{(1 + y_j) x_i}} \cdot \frac{x_i}{1+t} \sum_k \frac{x_k}{t (1 + y_i)} \right)$$
$$ = b_{ij} \int_{\Omega \setminus \Omega(u_i^{-1})} P(s, \omega) d\mu(\omega) \ + \frac{C}{1+t} [(b_{ij} (1
+ y_j) - (t - x_i - x_j) b_{ij} - b_{ij} ]$$
Since $\mu(\Omega(u_i^{-1})) = x_i/(1+t),$ the second term is simply the $b_{ij}$ times the integral of
$P(s, \cdot)$ over $\Omega(u_i^{-1})$. We are done with the first part of the proof. 

\

2. Next we show that $\phi(u_i s) = a_i \phi(s)$ if $s$ is a reduced word not beginning with $u_i^{-1}$ (including
the possibility $s = 1$). To begin with,
$$\phi(u_i s) = \ <\pi(s) {\sf 1}, \pi(u_i^{-1}) {\sf 1} > \ = \int_\Omega P(s, \omega) P(u_i^{-1}, \omega)
d\mu(\omega)$$
$$ = \int_{\Omega \setminus \Omega(u_i^{-1})} + \sum_{j \not = i} \int_{\Omega(u_i^{-1} u_j)} + \sum_j
\int_{\Omega(u_i^{-1} u_j^{-1})} P(s, \omega) P(u_i^{-1}, \omega) d\mu(\omega) \ .$$
Use Lemma \ref{constant}(b) to find $C$ such that $P(s, \omega) = C$ for $\omega$ in  $\Omega(u_i^{-1})$.
Looking up values for $P(u_i^{-1}, \omega)$ in the various cases, we obtain
$$\phi(u_i s) = \sqrt{\frac{x_i}{t(1+y_i)}} \int_{\Omega \setminus \Omega(u_i^{-1})} P(s, \omega) d\mu(\omega) \
+$$
$$ C \left( - \sqrt{\frac{1+y_i}{t x_i}} \sum_{j \not= i} \hat{\mu}(u_i^{-1} u_j) + \sqrt{\frac{t(1+y_i)}{x_i}}
\sum_j \hat{\mu}(u_i^{-1} u_j^{-1}) \right) .$$
Evaluating the $\hat{\mu}$'s gives
$$\phi(u_i s) = a_i \int_{\Omega \setminus \Omega(u_i^{-1})} P(s, \omega) d\mu(\omega) \ +$$
$$C \left( - \sqrt{\frac{1+y_i}{t x_i}}  \cdot \frac{x_i}{1+t}  \sum_{j \not= i} \frac{x_j}{1+y_i} +
\sqrt{\frac{t(1+y_i)}{x_i}} \cdot \frac{x_i}{1+t} \sum_j \frac{x_j}{(1 + y_i) t} \right)$$
$$= a_i \int_{\Omega \setminus \Omega(u_i^{-1})} P(s, \omega) d\mu(\omega) \ + \frac{C a_i}{1 + t} (-y_i + t) $$
$$= a_i \left(\int_{\Omega \setminus \Omega(u_i^{-1})} + \int_{\Omega(u_i^{-1})} P(s,\omega) d\mu(\omega) \right)
\ ,$$
and we are done with the second part of the proof.

\

3. Finally we show that $\phi(s u_i) = a_i \phi(s)$ if $s$ is either 1 or a reduced word ending with $u_j$ for some
$j$. Let $\Omega(-) = \Omega(u_1^{-1}) \cup \Omega(u_2^{-1}) \cup \ldots \cup \Omega(u_n^{-1}).$ Then
$$\phi(s u_i) = \int_\Omega P(u_i, \omega) P(s^{-1}, \omega) d \mu(\omega) = $$
$$ = \int_{\Omega(u_i)} + \sum_{k \not= i} \int_{\Omega(u_k)} + \int_{\Omega(-)} P(u_i, \omega)
P(s^{-1}, \omega) d \mu(\omega).$$
By Lemma \ref{constant} (a), we have a $C$ such that $P(s^{-1}, \omega) = C$ in the first two integrals. Thus
$$\phi(s u_i) = C \left( \sqrt{\frac{t (1 + y_i)}{x_i}} \hat{\mu}(u_i) - \sqrt{\frac{t x_i}{1 + y_i}} \sum_{k \not=
i} \hat{\mu}(u_k) \right) + a_i \int_{\Omega(-)} P(s^{-1}, \omega) d \mu(\omega)$$
$$= C \left( \sqrt{\frac{t (1 + y_i)}{x_i}} \cdot \frac{x_i}{t (1 + t)} - \sqrt{\frac{t x_i}{1 + y_i}} \sum_{k \not = i}
\frac{x_k}{t (1 + t)} \right) + a_i \int_{\Omega(-)} P(s^{-1}, \omega) d \mu(\omega)  $$
$$= \frac{C a_i}{1+t} ((1+y_i) - y_i) + a_i \int_{\Omega(-)} P(s^{-1}, \omega) d \mu(\omega)$$
$$= a_i \left( \int_{\Omega(u_1) \cup \ldots \cup \Omega(u_n)} + \int_{\Omega(-)} P(s^{-1}, \omega) d \mu(\omega)
\right) = a_i \phi(s^{-1}) = a_i \phi(s)$$
since $\mu(\Omega(u_1) \cup \ldots \cup \Omega(u_n)) = \sum_k \beta(u_k) = 1/(1+t).$

As we observed in Lemma \ref{algprop}, parts 1, 2, and 3 above suffice to establish the asserted formula for
$\phi$. \end{proof}

\

We remark that the $a$'s and $b$'s above are related as in Section \ref{pure}. Set $b_{ii} =1$ for each $i$,
and let $B$ be the $n \times n$ matrix $(b_{ij})$. Then $B$ is positive because $\phi$ is positive definite and $b_{ij}
= \phi(u_i^{-1}u_j).$ Since
$$\frac{c_i}{\lambda} = \sqrt{\frac{x_i (1 + y_i)}{t}},$$
the $i^{\mbox{th}}$ entry of $B \vec{c}/\lambda$ is
$$\sqrt{\frac{x_i (1 + y_i)}{t}} + \sum_{j \not= i} b_{ij} \sqrt{\frac{x_j (1 + y_j)}{t}}
 = \sqrt{\frac{x_i (1 + y_i)}{t}} - \sqrt{\frac{x_i}{t (1 + y_i)}} \sum_{j \not= i} x_j$$
$$ = \sqrt{\frac{x_i}{t (1 + y_i)}} (1 + y_i - y_i) = a_i \ ,$$
and 
$$B \vec{c} \cdot \vec{c} = \lambda \sum_i a_i c_i = \lambda \sum_i c_i \phi(u_i) = \lambda^2$$
by Proposition \ref{rightsum}. We thus have the ingredients for one of the states examined in Section
\ref{pure}, and our $\phi$ is $\phi_{B, \vec{c}}$ for the $B$ we have described.

\

We can now record (most of) the main result of this paper; see also Proposition \ref{onlyev} below.

\

\begin{thm}\label{main} Given positive numbers $c_1, c_2, \ldots, c_n$ and $\lambda$ satisfying
$$\lambda^2 < \sum_i c_i^2 \ \ \mbox{and} \ \ c_j^2 - \sum_{i \not= j} c_i^2  < \lambda^2  \ \ \mbox{for} \ \ j =
1, \ldots, n \ ,$$ there are unique positive numbers $x_1, x_2, \ldots, x_n$ such that
$$\frac{x_j (1 + \sum_{i \not= j} x_i)}{\sum_i x_i} = \frac{c_j^2}{\lambda^2}$$
for each $j$. The function $\phi$ defined on $G$ by
$$\phi(s) = (- \sum_j x_j)^{\gamma(s)} \prod_i \left(\frac{x_i}{(1 + \sum_{j \not= i} x_j) \sum_j x_j} \right)
^{|s|_i/2} ,$$
where $|s|_i$ and $\gamma(s)$ are respectively the number of occurrences in $s$ of $u_i^{\pm1}$ and
$u_j^{-1}u_k$ for $k \not = j$, is a reduced pure $\lambda$-eigenstate for $\sum_i c_i u_i$.
\end{thm}

\

\begin{proof}  See Section \ref{mapS} below for the existence and uniqueness of the $x_i$'s. The
formula for $\phi$ follows from Proposition \ref{stateformula} above and the observation that the $a$'s and $b$'s
there satisfy $b_{ij} = - t a_i a_j$ for $i \not= j$. (Notice as well that $\gamma =  \sum_{j > k} \gamma_{jk}$.)
Proposition \ref{rightsum} implies that $\phi$ is a $\lambda$-eigenstate for $\sum_i c_i u_i$.
That $\phi$ is  pure follows from results in Section \ref{pure} and our identification of
$\phi$ as $\phi_{B, \vec{c}}$ above. Finally,  it follows from Theorem 2.7 in \cite{Spielberg} that $\phi$ is
reduced. \end{proof}

\

\section{Further observations}\label{further}

\

Somewhat surprisingly, the constant function ${\sf 1}$ turns out not to be cyclic for the boundary representation
$\pi$ described above. Let $H_1$ be the closed linear span of $\pi(G) {\sf 1}$, so  by
Proposition \ref{stateformula}, the restriction $\pi_1$ of $\pi$ to $H_1$ is one of the irreducible representations
described in Section \ref{pure} above. After some preparation, we will show that the restriction $\pi_2$ of $\pi$ to
the orthogonal complement $H_2$ of $H_1$ is unitarily equivalent to $\pi_1 \circ \sigma,$ where $\sigma$ is the
automorphism of $G$ that sends each $u_i$ to $u_i^{-1}$, and that the latter is unitarily inequivalent to $\pi_1$. This
behavior distinguishes $\pi_1$ sharply from the representations studied by Kuhn and Steger in \cite{KuhnSteger}. It
seems likely that $\pi_1$ belongs to the ``odd'' category in the tripartite classification of irreducible representations
proposed and conjectured to be exhaustive in \cite{KuhnSteger2}, but this is a matter that will have to be pursued
elsewhere.

Write $\Omega(\pm) = \bigcup_{j = 1}^n
\Omega(u_j^{\pm1})$ as above, and consider $$h = \sqrt{t} \ \chi_{\Omega(+)} - \frac{1}{\sqrt{t}} \
\chi_{\Omega(-)},$$ where the $\chi$'s are indicator functions. This is a unit vector in $L^2(\Omega, \mu)$ because
the measures of
$\Omega(+)$ and $\Omega(-)$ are respectively $1/(1+t)$ and $t/(1+t)$. We will see below that $H_2$ is the closed
$\pi(G)$-invariant subspace generated by $h$.

\

\begin{lem}\label{valP} For $v$ in $V$ and for $s$ in $G$ not ending in $v^{-1}$, we have
$$P(v,s^{-1} \omega) = \frac{\phi(sv)}{\phi(s)}$$
either for all $\omega$ in $\Omega(+)$ or for all $\omega$ in $\Omega(-)$. 
\end{lem}

\

\begin{proof} When $s = 1$, this is immediate from the descriptions in the vicinity of Proposition
\ref{rightsum} above for
$P(u_i^{\pm1}, \cdot).$ Suppose $s$ begins with $u_j^{\pm1}$, so there is no canceling in $s^{-1} \omega$ for
$\omega$ in $\Omega(\mp).$ If $s$ doesn't end in $u_i$, then
$$P(u_i^{-1}, s^{-1} \omega) = a_i = \frac{\phi(su_i^{-1})}{\phi(s)}$$
for all $\omega$ in $\Omega(\mp)$. This takes care of the case $v = u_i^{-1}$. If $s$ ends in $u_k$ for some $k$,
then 
$$P(u_i, s^{-1} \omega) = a_i = \frac{\phi(su_i)}{\phi(s)}$$
for all $\omega$ in $\Omega(\mp)$, while if $s$ ends in $u_k^{-1}$ for some $k$ different from $i$, we have
$$P(u_i, s^{-1} \omega) = -ta_i = \frac{b_{ki}}{a_k} = \frac{\phi(su_i)}{\phi(s)}$$
for all $\omega$ in $\Omega(\mp)$, which finishes the case $v = u_i$. \end{proof}

\

\begin{lem}\label{perp} The vector $h$ is orthogonal to $\pi(G) {\sf 1}$.
\end{lem}

\

\begin{proof} We must show that
$$(*) \ \ \ \ \ \ \int_{\Omega(-)} P(s, \omega) \, d\mu(\omega) = t \int_{\Omega(+)} P(s, \omega) \, d\mu(\omega)$$
for all $s$ in $G$. This is immediate when $s = 1,$ since $\mu(\Omega(u_j^{-1})) = t \mu(\Omega(u_j))$ for each
$j$. Suppose that (*) holds for some given $s$ in $G$ and that $v$ in $V$ is such that $sv$ is reduced. We have
$$\int_{\Omega(-)} \ + \int_{\Omega(+)} P(sv, \omega) \, d\mu(\omega) = \phi(sv) = \frac{\phi(sv)}{\phi(s)} \left(
\int_{\Omega(-)} \ + \int_{\Omega(+)} P(s, \omega) \, d\mu(\omega) \right) $$
because $\phi(s)$ is the integral of $P(s,\cdot)$ over $\Omega$.
It follows from Lemma \ref{valP} and $P(sv,\omega) = P(s,\omega) P(v, s^{-1} \omega)$ that either the
$\Omega(-)$ summands above are equal, or the $\Omega(+)$ summands are. Hence both are equal, and (*) for
$sv$ follows from (*) for $s$ by multiplying by $\phi(sv)/\phi(s).$ \end{proof}

\

We will use the following lemma (which likely holds in much greater generality) in proving the inequivalence of
$\pi_1$ and $\pi_2.$

\

\begin{lem}\label{invertible} Let $Y$ be a linear combination of 1 and $u_1, \ldots, u_n$ for which there exists $f$ in
$\ell^2(G)$ such that $Y \ast f = \delta_1.$ Then $Y$ is invertible in the left regular representation.
\end{lem}

\

\begin{proof} Write $Y = d_0 + \sum_i d_i u_i,$ where without loss of generality the $d$'s are all nonzero.
Since $d_0 f(1) + \sum_i d_i f(u_i^{-1}) = 1,$ at least one of $f(1), f(u_1^{-1}), \ldots, f(u_n^{-1})$ must be different
from zero. Suppose that $f(1) \not= 0$. For each index $j$, let $g_j$ be $f$ times the indicator function of $S_j^-$ (the
set of reduced words ending in $u_j^{-1}$). One checks readily that $Y \ast g_j = d_j f(u_j^{-1}) \delta_1.$ Thus, $g_j =
d_j f(u_j^{-1}) f$ \cite{Linnell}. Since $g_j(1) = 0 \not= f(1),$ this makes $g_j = 0;$ that is, $f$ vanishes on each $S_j^-$.
Now fix distinct indices $j$ and $k$, and let $r(s) = f(s u_k)$ if $s$ in $S_j^-$ and 0 otherwise. As with the $g_j$'s, we
have $Y \ast r = d_j r(u_j^{-1}) \delta_1,$ so $r$ must be a multiple of $f$, but $r(1) = 0,$ so $r = 0.$ Continuing in
this fashion, we see that $f$ must be supported on $G^+$, and it then follows by equating coefficients that $f(1) =
1/d_0, f(u_j) = -d_j/d_0^2,$ and in general
$$f(u_{i_1} u_{i_2} \ldots u_{i_m}) = (-1)^m\frac{d_{i_1} d_{i_2} \ldots d_{i_m}}{d_0^{m+1}} \ .$$
This makes
$$|d_0|^2 ||f||_2^2 = \sum_{m=0}^\infty \left( \frac{|d_1|^2 + \ldots + |d_n|^2}{|d_0|^2} \right) ^m,$$
so $|d_1|^2 + \ldots + |d_n|^2 < |d_0|^2,$ so $Y$ is invertible by Proposition \ref{redspectrum} above.

If $f(1) = 0,$ then $f(u_j^{-1}) \not= 0$ for some $j$. We have
$$(d_j + d_0 u_j^{-1} + \sum_{i \not= j,0} d_i u_i u_j^{-1}) \ast (u_j \ast f) = \delta_1,$$
and $(u_j \ast f)(1) = f(u_j^{-1}),$ so we may apply the previous case to the free generators $u_j^{-1}, u_i u_j^{-1} (i
\not= j).$ \end{proof}

\

\begin{prop}\label{twopieces} \ \ \ (a) The subspace $H_2$ is the closed linear span of $\pi(G) h$. 

(b) The representation $\pi_2$ is unitarily equivalent to $\pi_1 \circ \sigma,$ where $\sigma$ is the automorphism
of $G$ sending each $u_i$ to $u_i^{-1}$. 

(c) The representations $\pi_1$ and $\pi_2$ are unitarily inequivalent.
\end{prop}

\

\begin{proof} Consider the symmetry $T: \Omega \rightarrow \Omega$ taking each string in $\Omega$ to
the string obtained by inverting each symbol. It is immediate that $T(\Omega(s)) = \Omega(\sigma(s)).$
Since
$$\frac{\hat{\mu}(\sigma(v_1 v_2 \ldots v_k))}{\hat{\mu}(v_1 v_2 \ldots v_k)} = \frac{\beta(v_1^{-1})
\alpha(v_1^{-1},v_2^{-1}) \ldots \alpha(v_{k-1}^{-1},v_k^{-1})}{\beta(v_1) \alpha(v_1,v_2) \ldots
\alpha(v_{k-1},v_k)} = \frac{\beta(v_1^{-1})}{\beta(v_1)}$$
for a nonempty reduced word $v_1 v_2 \ldots v_k$ with the $v$'s in $V$, it follows that
$$\frac{d \mu (T \omega)}{ d \mu(\omega)} = \left\{\begin{array}{ll} t & \ \ \omega \in \Omega(+) \\
1/t & \ \ \omega \in \Omega(-) \end{array} \right. \ ,$$
so $h^2 = d(\mu \circ T)/d \mu.$ The operator $W$ on $L^2(\Omega, \mu)$ defined by 
$$(W f) (\omega) = h(\omega) f(T \omega)$$
is unitary and takes $\sf 1$ to $h$. We claim that $W$ intertwines $\pi$ and $\pi \circ \sigma$. For $f$ in
$L^2(\Omega, \mu),$ we have
$$(W \pi(u_i) f) (\omega) = h(\omega) P(u_i, T \omega) f(u_i^{-1} T \omega),$$
while
$$(\pi(u_i^{-1}) W f) (\omega) = h(u_i \omega) P(u_i^{-1}, \omega) f(u_i^{-1} T \omega) \ .$$
It is readily checked that $h(\omega) P(u_i, T \omega)$ and $h(u_i \omega) P(u_i^{-1}, \omega)$ are the same for
all $\omega$, namely
$$\left\{ \begin{array}{ll} -1/(\sqrt{t} \ a_i) & \ \ \ \omega \in \Omega(u_i^{-1}) \\
\sqrt{t} \ a_i & \ \ \ \omega \in \Omega \setminus \Omega(u_i^{-1}) \end{array} \right. \ \ .$$
(One must check two subcases in each of the two cases above.) The intertwining of each $\pi(u_i)$ with its inverse,
and thus of $\pi$ with $\pi \circ \sigma,$ now follows. 

Now let $K$ be the closed $\pi(G)$-invariant subspace of $L^2(\Omega,\mu)$ generated by ${\sf 1}$ and $h$. Parts (a)
and (b) will both follow once we show that $K = L^2(\Omega,\mu)$. Since ${\sf 1}$ and $h$ are linearly independent
linear combinations of $\chi_{\Omega(+)}$ and $\chi_{\Omega(-)}$, these two indicator functions must both belong
to $K$. We have
$$\pi(u_i) {\sf 1} = a_i \chi_{\Omega(-)} + \frac{1}{a_i} \chi_{\Omega(u_i)} - t a_i \sum_{j \not= i}
\chi_{\Omega(u_j)},$$
so
$$\frac{1}{a_i} \chi_{\Omega(u_i)} - t a_i \sum_{j \not= i} \chi_{\Omega(u_j)} \in K$$
for each $i$. It will follow that each $\chi_{\Omega(u_i)} \in K$ once we show that the matrix
$$\left( \begin{array}{cccccc} 
1/a_1 & -t a_1 & -t a_1 & \ldots & -t a_1 & -t a_1 \\
-t a_2 & 1/a_2 & -t a_2 & \ldots & -t a_2 & -t a_2 \\
&& . &&\\
&& . &&\\
&& . &&\\
-t a_n & -t a_n & -t a_n & \ldots & -t a_n & 1/a_n \end{array} \right) $$
is invertible. Multiplying  row $i$ by $a_i (1 + y_i)$ for each $i$ and using $ta_i^2 = x_i/(1+y_i)$ turns
this into
$$\left( \begin{array}{cccccc} 
1+y_1 & -x_1 &-x_1 & \ldots &-x_1 &-x_1 \\-x_2 & 1+y_2 &-x_2 & \ldots &-x_2 &-x_2 \\
&& . &&\\
&& . &&\\
&& . &&\\-x_n &-x_n &-x_n & \ldots &-x_n & 1+y_n \end{array} \right) \ .$$
whose determinant by Lemma \ref{detM} in the appendix is 
$$(1+t)^n - (1+t)^{n-1} \sum_j x_j = (1+t)^{n-1} \ .$$
We conclude that $\chi_{\Omega(u_i)} \in K$ for each $i$. The unitary $W$ preserves $K$ and sends 
$\chi_{\Omega(u_i)}$ to $t^{-1/2} \chi_{\Omega(u_i^{-1})}$, so $\chi_{\Omega(u_i^{-1})} \in K$ as well. Suppose we
have shown for some $k \geq 1$ that $\chi_{\Omega(s)} \in K$ for all reduced words $s$ of length $k$. For
such a word $s$ not beginning with $u_i^{-1},$ we have
$$\pi(u_i) \chi_{\Omega(s)} = P(u_i, \cdot) \chi_{\Omega(u_i s)} = \frac{1}{a_i} \chi_{\Omega(u_i s)},$$
so $\chi_{\Omega(u_i s)} \in K.$ Using the unitary $W$, it follows that $\chi_{\Omega(u_i^{-1} s)} \in K$ as well if
$s$ doesn't begin with $u_i$. We conclude that the indicator function of every cylinder set belongs to $K$,
showing that $K = L^2(\Omega, \mu).$ This takes care of (a) and (b).

To prove (c), it will suffice to show that $\sum_i c_i u_i^{-1}$ cannot have $\lambda$ as an eigenvalue in $\pi_1.$ We
will think of $\pi_1$ as acting on the Hilbert space
$$H = H_0 \oplus \bigoplus_{i = 1}^n \left(\ell^2 (S^-_i) \otimes H_i' \right) $$
constructed in Section \ref{pure} with $a_i$ and $b_{ij}$ as in Proposition \ref{stateformula}. Thus,  $H_0$ is the closed
linear span of $\pi_1(G^+) \Delta_1$, while $H_i'$ is  $H_0 \ominus \pi_1(u_i) H_0,$ and so forth. Suppose now that
$\xi$ in $H$ satisfies
$$\sum_i c_i \pi_1(u_i)^* \xi = \lambda \xi.$$
We claim first that $\xi$ must be orthogonal to each subspace $\ell^2(S_j^-) \otimes H_j'$. Fix $j$, and $\eta$ in
$H_j'$. Define $f$ in $\ell^2(G)$ by 
$$f(s) = \left\{\begin{array}{ll} <\xi, \delta_s \otimes \eta> & s \in S_j^- \\
0 & \mbox{else} \end{array} \right . \ .$$
Notice that $\lambda <\xi, \delta_s \otimes \eta> \ = \sum_i c_i <\xi, \pi(u_i) (\delta_s \otimes \eta)>$ for every $s$
in $G$. We obtain
$$ \lambda f(s) - \sum_i c_i f(u_i s) = \left\{\begin{array}{ll} c_j <\xi, \eta> & s = u_j^{-1}\\
0 & \mbox{else} \end{array} \right . $$
by checking the cases $s \in S_j^- \setminus \{u_j^{-1}\}$ (use $\pi_1(u_j) \delta_s \otimes \eta = \delta_{u_j s}
\otimes \eta$), $s = u_j^{-1}$ (use $\pi_1(u_j) \delta_{u_j^{-1}} \otimes \eta = \eta$ and $f(1) = 0$), and $s \notin
S_j^-$ separately.  This means that  
$$(\lambda - \sum_i c_i u_i^{-1}) \ast (f \ast u_j) = c_j <\xi, \eta> \delta_1 .$$
The assumption is ambient that $\lambda - \sum c_i u_i$ (and hence its adjoint) is not invertible in the left regular
representation, so $<\xi, \eta> = 0$ by Lemma \ref{invertible}, and thus $f = 0$ by 1.4 in \cite{Linnell}. 

We have so far shown that $\xi \in H_0,$ and that $\xi \perp H_j'$ (that is, $\xi \in \pi_1(u_j) H_0$) for each $j$.
Take $i \not= j$, for instance $i = 1$ and $j=2$. Using part (c) of Lemma \ref{orthdec}, we obtain 
$$\xi = \ <\xi,\Delta_{u_i}> \Delta_{u_i} =  \ <\xi,\Delta_{u_j}> \Delta_{u_j}.$$
This forces $\xi = 0,$ because $\Delta_{u_i}$ cannot be a scalar multiple of $\Delta_{u_j}$. (The scalar in question
would have to have modulus 1 because these are unit vectors, and would have to be $b_{ij}$ to make the inner
product with $\Delta_{u_j}$ come out right. However, $|b_{ij}| < 1$ because $b_{ij}$ is given by the formula in
Proposition \ref{stateformula}.) \end{proof}

\

A reasonable guess about the irreducible representations of $G$ treated in Section \ref{bdryrepns} is that they are
classified up to unitary equivalence by the vector $\vec{c}/\lambda.$ We leave this unresolved for now except to
note that if one fixes $\vec{c}$ and changes $\lambda$, the new representation is inequivalent to the original one.
This is an immediate consequence of the following proposition.

\

\begin{prop}\label{onlyev} The only eigenvalue of $\sum c_i \pi_1(u_i)$ is $\lambda$. The corresponding
eigenspace is ${\mathbb C} \Delta_1.$
\end{prop}

\

\begin{proof} We use the notation of the proof of part (c) of the previous proposition. Suppose $\sum_i c_i \pi_1(u_i)
\xi = \nu \xi$ for some complex number $\nu$ and some nonzero $\xi$ in $H$. Fix an index $j$, pick $\eta$ in
$H_j'$ and consider $g$ in $\ell^2(G)$ defined by 
$$g(s) = \left\{\begin{array}{ll} <\xi, \delta_s \otimes \eta> & s \in S_j^- \\
0 & \mbox{else} \end{array} \right . \ .$$
One checks readily that 
$$(\sum_i c_i u_i -\nu) \ast g = c_j g(u_j^{-1}) \delta_1.$$
Since $\nu$ belongs to the reduced spectrum of $\sum_i c_i u_i,$ it follows from Lemma \ref{invertible} above that
$g(u_j^{-1}) = 0,$ so $g$ vanishes identically by \cite{Linnell}. We have shown that $\xi$ must belong to $H_0$.
By Lemma \ref{orthdec}, then, we have
$$\nu <\xi, \Delta_1> \ = \sum_i c_i <\xi, \pi_1(U_i)^* \Delta_1> \ = \sum_i c_i a_i <\xi, \Delta_1> = \lambda <\xi,
\Delta_1> \ ,$$ 
and similarly for $s$ in $G^+$ and any index $j$
$$\nu <\xi, \Delta_{u_j s}> \ = c_j <\xi, \Delta_s> + \phi(s) \left( \sum_{i \not= j} c_i b_{ij} \right) <\xi, \Delta_1>.$$
We can't have $<\xi, \Delta_1> \ = 0,$ because that would force $<\xi, \Delta_s> = \ 0$ for all $s$ in $G^+$ and hence
$\xi =0.$ Hence $\nu = \lambda$. We may assume that $<\xi, \Delta_1> \ = 1.$ The preceding formula becomes
$$\lambda <\xi, \Delta_{u_j s}> \ = c_j <\xi, \Delta_s> + \phi(s) (\lambda a_j - c_j),$$
whence it readily follows by induction on the length of $s$ that $<\xi, \Delta_s> = \phi(s)$ for every $s$ in $G^+$, and
thus that $\xi = \Delta_s.$ \end{proof}

\

We have so far avoided spectral values on the boundary of the reduced spectrum. The situation there is simpler and
more clear-cut than in the interior (as well as being qualitatively different in the sense of \cite{KuhnSteger2}). 

\

\begin{thm}\label{unique} Let $\lambda = (\sum_i c_i^2)^{1/2}$, where $c_1, \ldots, c_n$ are nonnegative and not all
zero. The function $\phi$ defined on $G$ by 
$$\phi(s) = \left\{\begin{array}{ll}
\prod_i
\left(\frac{c_i}{\lambda}\right)^{|s|_i} & \mbox{if} \ \ \gamma(s) = 0 \\ 0 & \mbox{else} \end{array} \right . \ ,$$
where we understand $0^0 = 1$, is the unique reduced $\lambda$-eigenstate for $\sum_i c_i u_i$. 
\end{thm}

\

\begin{proof} If only one of the coefficients is nonzero, this is Lemma 4.4 in \cite{Paschke}. Assume therefore that
at least two coefficients are nonzero. The argument from \cite{Paschke} for the case $c_1 = \ldots = c_n$ also works
here with just a few cosmetic changes. 

Let $T = \lambda^{-1} \sum_i c_i u_i,$ thought of as an operator on $\ell^2(G).$  Because $1$ belongs to the reduced
spectrum of $T,$ there is a state $f$ on the algebra of bounded operators on $\ell^2(G)$ such that $f((T^*-1)(T-1)) =
0.$ We will be done once we show that the restriction of $f$ to $G$ must coincide with $\phi$. Let $S^+$ be the set of
reduced words in $G$ beginning with some $u_j$, and let $S^- = G \setminus S^+$. Let $P$ and $Q$ be respectively the
orthogonal projections of $\ell^2(G)$ on $\ell^2(S^+)$ and $\ell^2(S^-)$, so $Q = 1- P$. Suppose we know that $f(P) = 1.$
For an $s$ in $G$ with $\gamma(s) > 0$ (that is, for $s$ not in $G^+(G^+)^{-1}$), we have $P(T^*)^m s T^m P = 0$ for
sufficiently large $m$, and hence $f(s) = f(P(T^*)^m s T^m P) = 0$. On the other hand, if $s \in G^+(G^+)^{-1}$ and
$su_j^{-1}$ is reduced, then 
$$f(su_j^{-1}) = f(su_j^{-1} T) = \frac{c_j}{\lambda} f(s)$$
because $\gamma(s u_j^{-1} u_i) = 0$ for $i \not= j$. Likewise 
$$f(u_j s) = \frac{c_j}{\lambda} f(s)$$
if $u_j s$ is reduced. It now follows easily that $f(s) = \phi(s)$ for $s$ in as well as outside of $G^+(G^+)^{-1}.$

We show now that $f(P)$ must be $1$. Suppose not, that is, suppose $f(Q) > 0.$ Consider the state $g$ defined on
bounded operators $X$ by $g(X) = f(QXQ)/f(Q),$ so $g(Q) = 1.$ We have $QTT^*Q = Q, \ QT^*Q = T^*Q,$ and $QTQ = QT.$
Because $T$ is in the left kernel of $f$, this makes $g((T - 1) (T^* - 1)) = 0.$ The same argument as in the previous
paragraph, mirror-imaged by the automorphism $\sigma$ of $G$ that sends each $u_i$ to $u_i^{-1}$, shows that
$g = \phi \circ \sigma.$ In particular, we have $g(u_i^{-1}u_j) = \lambda^{-2} c_i c_j$  for $i \not= j.$ Consider now
$f(u_i^{-1}u_j)$. We have already observed that 
$$f((Q - QT)(Q - T^*Q)) = f(Q) g((T - 1) (T^* - 1)) = 0,$$
so $f(Pu_i^{-1}u_jQ) = f(Pu_i^{-1}u_jT^*Q) = f(0) = 0,$ and likewise $f(Qu_i^{-1}u_jP) = 0.$ Since also $Pu_i^{-1}u_jP =
0,$ it follows that 
$$f(u_i^{-1}u_j) = f(Qu_i^{-1}u_jQ) = f(Q) g(u_i^{-1}u_j) = f(Q) \lambda^{-2} c_i c_j$$
for $i \not= j$. All of these quantities are nonnegative, and at least two are positive. This, however, contradicts
$f(T^*T) = f(1) = 1,$ because the latter forces the sum over unequal $i$ and $j$ of $c_i c_j f(u_i^{-1}u_j)$ to
vanish. \end{proof}

\

Rotating the generators gets the uniqueness assertion above (with appropriately modified state formula) for an
arbitrary nonzero coefficient vector $\vec{c}$ and complex $\lambda$ with $|\lambda| = |\vec{c}|$.  As for the inner
boundary of the spectral annulus, if $|c_j|$ is the maximum of the absolute coefficients, and if
$$|c_j|^2 > \sum_{i \not= j} |c_i|^2 \ \ \ \mbox{and} \ \ \ |\lambda|^2 = |c_j|^2 - \sum_{i \not= j} |c_i|^2,$$
then $\sum_i c_i u_i$ has a unique reduced $\lambda$-eigenstate because 
$$\lambda u_j^{-1} - \sum_{i \not= j} c_i u_j^{-1} u_i$$
has a unique reduced $c_j$-eigenstate.

\

\section{Appendix --- the map S}\label{mapS}

Let ${\mathbb R}^n_+$ denote the positive orthant of ${\mathbb R}^n$. In this section we study the map $S: {\mathbb R}^n_+ \rightarrow
{\mathbb R}^n_+$ defined by
$$S(x_1, x_2 , \ldots x_n) = \frac{1}{t} (x_1(1 + y_1), x_2(1+ y_2), \ldots, x_n(1 + y_n)),$$
where
$$t =  \sum_i x_i\ \ \mbox{and for each} \ j , \ y_j = \sum_{i \not= j} x_i \ .$$
We will use this notation  --- $t$ for the sum of the $x$'s, $y_j$ for $t - x_j$  --- 
throughout our discussion of $S$. Our goal is to show that $S$ is one-to-one on ${\mathbb R}^n_+$ and that the range
of $S$ is the open subset $D_n$ of ${\mathbb R}^n_+$ consisting of the points $(s_1, s_2, \ldots s_n)$ in
${\mathbb R}^n_+$ such that
$$\sum_i s_i > 1, \ \ \mbox{and for each} \ j , \ s_j < 1  \ + \sum_{i \not= j} s_i \ .$$
That $S({\mathbb R}^n_+) \subseteq D_n$ is elementary. If $s_i = x_i(1 + y_i)/t$ for each $i$, then
$$\sum_i s_i \ > \ \frac{1}{t} \sum_i x_i \ = \ 1,$$
and for each $j$,
$$t \left( \sum_{i \not= j} s_i + 1 - s_j \right) =  \sum_{i \not= j} x_i(1+t-x_i) \ + \ t \ - \ x_j(1+t-x_j) $$
$$ = (t - x_j) (1+ t) - \sum_{i \not= j}x_i^2 + t - x_j - t x_j + x_j^2 = 2 (t - x_j) + (t - x_j)^2 - \sum_{i
\not= j} x_i^2$$
$$= 2 \sum_{i \not= j} x_i + \left(\sum_{i \not= j} x_i \right)^2 - \sum_{i \not= j} x_i^2 \ > \ 0 \ .$$

The derivative of $S$ is easily calculated. For $i \not= j$, the $(i,j)$ entry is given by
$$S'_{i,j} = \frac{\partial}{\partial x_j} \frac{x_i (1 + y_i)}{t} = \frac{x_i^2 - x_i}{t^2},$$
while the diagonal entries are given by
$$S'_{i,i} = \frac{y_i^2 + y_i}{t^2} \ .$$ 
The following lemma will help show that $\det S'$ is positive on ${\mathbb R}^n_+$

\

\begin{lem}\label{detM} Let $M$ be an $n \times n$ matrix of the form
$$\left( \begin{array}{cccccc} 
r_1 & p_1 & p_1 & \ldots & p_1 & p_1 \\
p_2 & r_2 & p_2 & \ldots & p_2 & p_2 \\
&& . &&\\
&& . &&\\
&& . &&\\
p_n & p_n & p_n & \ldots & p_n & r_n \end{array} \right) . $$
Let $q_j = r_j - p_j$ for each $j$. The determinant of $M$ is
$$ \prod_i q_i + \sum_j p_j \prod_{i \not= j} q_i$$
\end{lem}

\

\begin{proof} Subtract the first column of $M$ from the other columns to obtain 
$$\det M = (r_1 q_2 q_3 q_4 \ldots q_n) + (p_2 q_1 q_3 q_4 \ldots q_n) + (p_3 q_2 q_1 q_4 \ldots q_n) + \ldots$$
$$+ (p_n q_2 q_3 q_4 \ldots q_{n-1} q_1),$$ 
then put $q_1 + p_1$ for $r_1$ in the first term. \end{proof}

\

\begin{lem}\label{JacS} \ \ \ $\det (S') > 0$ on ${\mathbb R}^n_+$ \ .
\end{lem}

\

\begin{proof} Notice that
$$y_j^2 + y_j - (x_j^2 - x_j) = (t- x_j)^2 - x_j^2 + t = t (1 + t - 2 x_j) = t (1 + y_j - x_j) \ .$$
Let $\theta_j = 1 + y_j - x_j$ . Apply the lemma above to $\det (t^2 S')$ and then divide by $t^{n-1}$ to write 
$$t^{n+1} \det(S') \equiv P(x_1, x_2, \ldots, x_n)$$
$$ = (t \theta_1 \theta_2 \theta_3  \ldots \theta_n) + (x_1^2 - x_1)(\theta_2 \theta_3 \ldots \theta_n) +
\theta_1 (x_2^2 - x_2) (\theta_3 \ldots \theta_n) + $$
$$\ldots + (\theta_1 \theta_2 \ldots \theta_{n-1}) (x_n^2 - x_n) .$$
Write the first term as 
$$\sum_j (x_j \theta_2 \ldots \theta_n)$$
and add the $j$th term of this sum to the $j$th of the remaining terms of $P$ to obtain 
$$P = \sum_j (x_j \theta_j + x_j^2 - x_j) \prod_{i \not= j} \theta_i \  = \ \sum_j x_j y_j \prod_{i \not=
j} \theta_i.$$ If all of the $\theta$'s are positive, we conclude immediately that $P > 0.$ Otherwise, since at
most one $x_i$ can exceed the sum of the others, all but one of the $\theta$'s must be positive.  Assume for
definiteness that $\theta_i > 0$ for $i \geq 2.$ Writing $y_1$ as $x_2 + \ldots + x_n$, we have
$$P = \sum_{j \geq 2} x_1 x_j \prod_{i \geq 2} \theta_i \ + \ \sum_{j \geq 2} x_j y_j \prod_{i \not=
j} \theta_i \ $$
$$ = \sum_{j \geq 2} x_j (x_1 \theta_j + y_j \theta_1) \prod_{i \not= 1,j} \theta_i \ .$$
Now 
$$x_1 \theta_j + y_j \theta_1 = x_1 + y_j + y_1 y_j - x_1 x_j \ ,$$
which is positive for $j \geq 2$ because $y_1 > x_j$ and $y_j > x_1.$ \end{proof}

\

Recall the set $D_n$ defined by strict linear inequalities at the beginning of this section.

\

\begin{prop}\label{rangeS} The map $S$ takes ${\mathbb R}^n_+$ onto $D_n.$
\end{prop}

\

\begin{proof} We have seen by now that $S$ is an open map of ${\mathbb R}^n_+$ into $D_n.$ It will suffice to
show that if $\vec{\sigma} = (\sigma_1,\sigma_2, \ldots, \sigma_n)$ is the limit of a sequence
$\{S(\vec{x}^{(m)})\}$, where
$\{\vec{x}^{(m)}\}_m$ is a sequence in ${\mathbb R}^n_+$ with no limit points in ${\mathbb R}^n_+$, then $\vec{\sigma} \notin
D_n$. Omit the superscript $m$ and write
$$s_j = \frac{x_j (1 + y_j)}{x_j + y_j}, \ s_j \rightarrow \sigma_j \ \ (j = 1, \ldots, n).$$
After passing to a subsequence, we may assume that either $x_j \rightarrow 0$ for some $j$, or $x_j \rightarrow
\infty$ for some $j$. Suppose that $x_j \rightarrow 0$ for some $j$. Notice that this forces $x_j \not= s_j$
(eventually along the sequence), else $x_j = s_j = 1.$ If $\sigma_j = 0,$ then $\vec{\sigma} \notin D_n$, and we are
done. If $\sigma_j > 0,$ we may write
$$y_j = \frac{x_j (s_j - 1)}{x_j - s_j}$$
and conclude that $y_j \rightarrow 0.$ This makes $x_i \rightarrow 0$ for every $i$, and thus $y_i \rightarrow
0$ for every $i$. Since
$$\frac{x_i}{x_i + y_i} + y_i \frac{x_i}{x_i + y_i} \rightarrow \sigma_i \, $$
we must have $x_i/(x_i + y_i)  \rightarrow \sigma_i$ for every $i$. Since $\sum_i x_i/(x_i + y_i) = 1,$ it
follows that $\sum_i \sigma_i = 1,$ so $\vec{\sigma} \notin D_n.$ The remaining possibility is that
$x_j \rightarrow \infty$ for some $j$. It follows that $y_j \rightarrow \sigma_j - 1.$ Since $y_i \rightarrow
\infty$ for $i \not= j,$ we also have 
$$x_i = \frac{s_i y_i}{1 + y_i - s_i} \rightarrow \sigma_i$$
for $i \not = j$. But $y_j = \sum_{i \not= j} x_i,$ so
$$\sigma_j - 1 = \sum_{i \not= j} \sigma_i \ ,$$
which rules $\vec{\sigma}$ out of $D_n$. \end{proof}

\

Showing that $S$ is injective on ${\mathbb R}^n_+$ is more troublesome than identifying its range. We begin with a simple
observation.

\begin{lem}\label{sumofentries} The sum of the entries of $S$ is a strictly increasing function of each of $x_i$.
\end{lem}

\

\begin{proof} The partial derivative with respect to $x_i$ of the sum of the entries of $S$ is
$$t^{-2} (y_i^2 + y_i + \sum_{j \not= i} (x_j^2 - x_j) ) = t^{-2} (y_i^2 + \sum_{j \not= i} x_j^2) > 0.  $$ \end{proof}

\

Consider now the situation that would arise if a point $\vec{s} = (s_1, \ldots, s_n)$ in $D_n$ were hit under
$S$ by two distinct points in ${\mathbb R}^n_+$, say $\vec{x}$ and $\vec{\xi}$. 

\

\begin{lem}\label{distinct} If two distinct points in ${\mathbb R}^n_+$ have the same image under $S$, then summing the
coordinates of these points yields two distinct numbers greater than 1.
\end{lem}

\

\begin{proof} Let $t$ (as usual) be the sum of the entries of $\vec{x}$, and let $\tau$ be the same for
$\xi$. Further let $\vec{b} = \vec{x}/t$ and $\vec{\beta} = \vec{\xi}/\tau.$  It follows easily that
$$ b_i + t (b_i - b_i^2) = s_i = \beta_i + \tau (\beta_i - \beta_i^2)$$ for each $i$, and further manipulation
shows that
$$(*) \ \ \ (\tau - t) (\beta_i - \beta_i^2) = (b_i - \beta_i)(1 + t (1 - b_i - \beta_i)) \ .$$ Since $\sum b_i = 1 = \sum
\beta_i,$ there is at most index $i$ such that $b_i + \beta_i > 1.$ If we had $\tau = t$, then by (*) $b_i$ and
$\beta_i$ would have to coincide for every $i$ except possibly the one for which $b_i + \beta_i > 1,$ which means
they must coincide for this $i$ as well. In order for $x$ and $\xi$ to be different, then, we must have $\tau
\not= t.$ Assume for definiteness that $\tau > t,$ so both sides of (*) are
greater than zero. Since the $b$'s and $\beta$'s both sum to 1, we must have $\beta_i > b_i$ for some $i$. For
this $i$, we must have $b_i + \beta_i > 1$ and $t > 1/(b_i + \beta_i - 1).$ Because $b_i + \beta_i < 2$, it follows that
$t$, and hence $\tau$, is greater than 1. \end{proof}

\

For $r > 0,$ let $\Delta_r = \{\vec{w} \in {\mathbb R}^n_+ : \sum w_i = r\}.$ As in the lemma and its proof, assume there
exist $\tau > t > 1$ such that $S(\Delta_\tau) \cap S(\Delta_t) \not= \emptyset.$ What then happens if we fix
$\tau$ and try to minimize $t$?

\

\begin{lem}\label{nomint} There is no minimal $t$ less than $\tau$ such that $S(\Delta_\tau) \cap S(\Delta_t)
\not= \emptyset.$
\end{lem}

\

\begin{proof} Suppose the pair ($\vec{x}, \vec{\xi}$) minimizes $\sum x_i \ ( = t)$ subject to $t < \tau,
S(\vec{x}) = S(\vec{\xi})$, and $\sum \xi_i = \tau$. Write $\vec{s} = S(\vec{x}) = S(\vec{\xi})$.  Let $\Phi$ be
a local inverse for $S$ such that $\Phi(\vec{s}) = \vec{x}$. (The existence of $\Phi$ follows from Lemma
\ref{JacS}.) Write $\vec{v}_1 = (1,1, \ldots, 1),$ and let $P$ be the hyperplane $\{\vec{v} \in {\mathbb R}^n : \vec{v} \cdot
\vec{v}_1 = 0 \}.$ For any $\vec{v}$ in $P$, the smooth real function $a \mapsto \Phi(S(\vec{\xi} + a \vec{v}))
\cdot \vec{v}_1$ is minimized by $a = 0$ and hence 
$$\Phi'(\vec{s}) S'(\vec{\xi}) \vec{v} \cdot \vec{v}_1 = 0 \ .$$ But of course $\Phi(S(\vec{x} + a \vec{v}))
\cdot \vec{v_1} \equiv t,$ so as well $$\Phi'(\vec{s}) S'(\vec{x}) \vec{v} \cdot \vec{v_1} = 0 \ .$$ 
It follows that the images of $P$ under $S'(\vec{\xi})$ and $S'(\vec{x})$ must coincide. (In other words, not
surprisingly, the hypersurfaces $S(\Delta_t)$ and $S(\Delta_\tau)$ must be tangent to one another at $\vec{s}$.)
The observation at the beginning of the proof of Lemma \ref{JacS} shows that $S'(\vec{x})$ acts on $P$ by
multiplying the $i^{\mbox{th}}$ entry of a vector in $P$ by $t(1 + t - 2 x_i)$. It follows that if $\vec{w}$ is
orthogonal to $S'(\vec{x})(P),$ then $w_i (1 + t - 2 x_i)$ is the same for all $i$. Describing the normal to
$S'(\vec{\xi})(P)$ in similar fashion, we deduce from $S'(\vec{x})(P) = S'(\vec{\xi})(P)$ that there is a real number
$\alpha$ such that
$$ (*) \ \ \ 1 + t - 2 x_i = \alpha ( 1 + \tau - 2 \xi_i)$$
for every $i$.  Summing on $i$, we obtain $n + (n-2) t = \alpha (n + (n-2) \tau).$ In case $n=2,$ this makes $\alpha =
1$ and (*) simply means that $x_1 - x_2 = \xi_1 - \xi_2.$ If we add this to the inequality $x_1 + x_2 < \xi_1 + \xi_2,$
we obtain $x_1 < \xi_1.$ Likewise, $x_2 < \xi_2.$ For $n > 2,$ we have $\alpha < 1,$ and the linear system (*) has the
form  $$A(\vec{x} - \alpha \vec{\xi}) = (\alpha - 1) \vec{v}_1,$$
where $A$ is easily seen to be invertible with $A^{-1} v_1 = (n-2)^{-1} v_1.$ Thus 
$$x_i = \alpha \xi_i + \frac{\alpha - 1}{n-2},$$
and hence in this case as well we have $x_i < \xi_i$ for each $i$. By Lemma \ref{sumofentries}, this contradicts
$S(\vec{x}) = S(\vec{\xi}).$ \end{proof}

\

What we have just shown will enable us to restrict attention to the behavior of $S$ near points in the boundary of
${\mathbb R}^n_+$ with all but one coordinate equal to zero. The next lemma treats that case in detail.

\

\begin{lem}\label{singleton} If $\vec{s} \in D_n$ and $\sum_{i \not= j} s_i < 1$ for some $j$, then $S^{-1}
(\{\vec{s}\})$ is a singleton.
\end{lem}

\

\begin{proof} Assume that $\sum_{i > 1} s_i < 1.$ Let $\vec{x}$ in ${\mathbb R}^n_+$ be such that $S(\vec{x}) = \vec{s}$,
and write $\vec{b} = \vec{x}/t.$ In light of Lemma \ref{distinct}, we may assume that $t > 1.$ As in the proof of
that lemma, we have
$$t b_i^2 - (1+t) b_i + s_i = 0 \ ,$$
so
$$b_i = \frac{1}{2t} (1 + t \pm \sqrt{(1+t)^2 - 4 t s_i})$$
for every $i$. For $i > 1,$ the sign choice must be $-$ rather than $+$, because $s_i < 1$ for such $i$ and the
choice of $+$ would make $b_i > 1.$ Let $Q_i$ be the
quadratic polynomial defined by
$$Q_i(r) = (1+r)^2 - 4 r s_i = r^2 + 2 (1 - 2 s_i) r + 1.$$
Notice that $Q_i$ has no real zeros when $i>1.$ Let $r_*$ be either 1 if $Q_1$ has no real zeros, or the larger
real zero of $Q_1$. That is,
$$r_* = \left\{ \begin{array}{cl}1 & \ \ \mbox{if} \ s_1 < 1\\ 2s_1 - 1 + 2 \sqrt{s_1^2 - s_1} & \ \ \mbox{if}
\ s_1 \geq 1 \end{array} \right. \ \ .$$
In either case, $t  \geq r_*$. (Since $s_1 > 1$ makes $Q_1(1) < 0$, the smaller zero of $Q_1$ is less than 1 in
this case.) Define functions $f_+$ and $f_-$ on the interval $[r_*, \infty)$ by 
$$f_{\pm}(r) = n + (n-2)r \pm Q_1(r)^{1/2} - \sum_{i>1} Q_i(r)^{1/2}.$$
Because the $b$'s sum to 1, we must have either $f_+(t) = 0$ or $f_-(t) = 0.$ For $r > r_*,$ the derivative of
$f_-$ is
$$f_-'(r) = n-2 - \frac{r + 1 - 2 s_1}{Q_1(r)^{1/2}} - \sum_{i > 1} \frac{r + 1 - 2 s_i}{Q_i(r)^{1/2}} \ .$$
Because $r > \max\{1,2 s_1 - 1\}, Q_i(r) < (r + 1)^2,$ and $\sum_{i>1} s_i < 1,$ we have
$$f_-'(r) < n - 2 - \sum_{i > 1} \frac{r + 1 - 2 s_i}{r + 1} = -1 + \frac{2}{r + 1} \sum_{i > 1} s_i < 0.$$
Thus $f_-$ is strictly decreasing on $[r_*, \infty).$ If $s_1 \leq 1,$ the sign choice $+$ is ruled out
for $b_1$ as well as for the other $b$'s, so we must have $f_-(t) = 0$ and the lemma is proved in this case.
Assume, then, that $s_1 > 1.$ We have 
$$f_+''(r) = 4 \frac{(s_1 - s_1^2)}{Q_1(r)^{3/2}} - 4 \sum_{i>1} \frac{s_i - s_i^2}{Q_i(r)^{3/2}},$$
all terms of which are negative. Since
$$\lim_{r \rightarrow \infty} f_+'(r) = n - 2 + 1 - (n-1) = 0$$
and $f_+'$ is strictly decreasing on $(r_*, \infty)$, it must be that $f_+$ is strictly increasing on $[r_*,
\infty)$. Since $f_+(r_*) = f_-(r_*)$ and the two functions in question are strictly monotone in opposite
senses, the condition that one or another vanish at $t$ determines $t$ uniquely. As we saw in the proof of Lemma
\ref{distinct}, this in turn determines $\vec{x}$ uniquely. \end{proof}

\

\begin{prop}\label{injective} The map S is injective on ${\mathbb R}^n_+.$
\end{prop}

\

\begin{proof} We proceed by induction on $n$. The case $n = 2$ is straightforward; one checks easily
that the map from $D_2$ to ${\mathbb R}^2_+$ defined by
$$(s_1,s_2) \mapsto \left(\frac{s_1 + s_2 - 1}{1 + s_2 - s_1},\frac{s_1 + s_2 - 1}{1 + s_1 - s_2} \right)$$
gives the identity map on ${\mathbb R}^2_+$ when preceded by $S$. Suppose now that $n > 2$ and that the assertion is
true in all dimensions less than $n$. Suppose that for some $\tau$, the set 
$$W = \{ (\vec{x},\vec{\xi}) \in {\mathbb R}^n_+ \times {\mathbb R}^n_+ : \sum_i x_i < \sum_i \xi_i = \tau, S(\vec{x}) =
S(\vec{\xi}) \}$$ is nonempty. Let $t_*$ be the infimum of the numbers $\sum_i x_i$ that come from pairs in $W$.
A sequence in $W$ along which the sum of the coordinates in the $\vec{x}-$slot tends to $t_*$ must by Lemma
\ref{nomint} have a subsequence converging to a pair $(\vec{x}_*, \vec{\xi}_*)$, where at least one coordinate in
each of
$\vec{x}_*$ and $\vec{\xi}_*$ vanishes (and because $S(\vec{x}) = S(\vec{\xi})$ along the sequence, the zeros
occur in same-indexed coordinates). If the number of nonzero coordinates in $\vec{x}_*$ (and $\vec{\xi}_*)$
is at least two, then the induction hypothesis is contradicted. If this number is one, then Lemma \ref{singleton} is
contradicted. There must of course be at least one nonzero coordinate because the sum of the coordinates of
$\vec{\xi}_*$ is $\tau.$  \end{proof}

\

We remark that the pleasing appearance of the formula for $S^{-1}$ in the case $n = 2$ does not at all reflect
what happens for larger $n$. For instance, as may easily be checked,
$$S^{-1} (1, 2, 3) = \left( 1, \frac{5 + \sqrt{73}}{6}, \frac{7 + \sqrt{73}}{2} \right) \ .$$

\

\

\end{document}